\documentclass{article}
\usepackage[utf8]{inputenc}
\usepackage{amsmath}
\usepackage{graphicx}
\usepackage{adjustbox}
\usepackage{url}
\usepackage{hhline}
\usepackage{float}
\restylefloat{table}

\title{Renormalization and blow-up for the 3D Euler equations}
\author{Jacob Price, Panos Stinis \\ 
University of Washington, Pacific Northwest National Laboratory}

\begin{document}
\maketitle
\begin{abstract}
In recent work we have developed a renormalization framework for stabilizing reduced order models for time-dependent partial differential equations. We have applied this framework to the open problem of finite-time singularity formation (blow-up) for the 3D Euler equations of incompressible fluid flow. To the best of our knowledge this is the first time-dependent perturbative renormalization approach for 3D Euler which includes all the complex effects present in the Euler dynamics. For the Taylor-Green initial condition, the renormalized coefficients in the reduced order models decay algebraically with time and resolution. The renormalized reduced models are stable and we evolve them for long times. Our results for the behavior of the solutions are consistent with the formation of a finite-time singularity. 
\end{abstract}

\section{Introduction}

The behavior of solutions of the 3D Euler equations for incompressible fluid flow is one of the most challenging problems in the analysis of partial differential equations (PDEs) and scientific computing. The common root of these difficulties lies in the complexity of the dynamics implied by the equations. On the theoretical side, it is unknown whether general smooth initial conditions give rise to a finite time singularity. On the numerical side, computations that begin with smooth initial conditions quickly give rise to a degree of complexity (turbulence) that exhausts the available computational power. Despite the accumulated knowledge on the theoretical and numerical fronts \cite{marchioro1994euler,constantin1994geometric,moffatt,C94,majda,deng2006improved,li2008blow,constantin2008euler,agafontsev2015development,tao2016blow,constantin2017analysis,isett2017onsager,elgindi2018finite,shu2005euler,hou2006dynamic,bustamante2012euler,grafke2013lagrangian,orlandi2012vortex,luo2014potentially,ayala2017extreme}, there is still no conclusive evidence.

A limiting factor in the numerical exploration of this problem is the necessary system size for fully resolved computations. Consequently, the development of reduced order models for Euler's equations would greatly help in this effort. A successful reduced order model for Euler's equations would describe a finite subset of the variables involved in a full scale simulation. For example, if the solution is represented in terms of Fourier modes, one could consider a finite subset of those modes corresponding to long wavelength behavior. Ideally, the reduced order models for these variables would evolve as the same variables would evolve should the full system have been simulated.

One mathematical framework for constructing reduced order models is the \emph{Mori-Zwanzig formalism} (MZ). Originally developed in the context of statistical mechanics \cite{zwanzig1961memory}, the formalism has been modernized as a mathematical tool \cite{chorin2002optimal,chorin2007problem}. This formalism allows one to express the evolution of a subset of variables (the \emph{resolved} variables) in terms of a Markov term, a noise term, and a memory integral. Based on various approximations, this framework has led to successful reduced order models for a host of systems (see e.g \cite{chorin2002optimal, bernstein2007optimal,stinis2012numerical,li2015memory,lei2016,parish2017}). 

Except for special cases, it is difficult to guarantee that the reduced models will remain stable. We have developed a time-dependent version of the \emph{renormalization} concept from physics \cite{goldenfeld1992,georgi1993} to aid with the stabilization of reduced models. This approach, in which we attach time-dependent coefficients to the terms in the reduced order model, has led to some success in stabilizing the reduced models \cite{stinis2013renormalized,stinis2015renormalized,price2017renormalized}. 

The MZ formalism has been previously used to develop reduced order models for Euler's equations \cite{hald2007optimal,stinis2007higher,stinis2013renormalized}. Those models were based on approximations of the memory term which assume various degrees of ``long memory'' i.e., the assumption that the unresolved variables evolve on timescales that are comparable to the resolved variables. Such an assumption is eminently plausible for the Euler equations (and high-Reynolds number fluid flows in general), given the vast range of active scales present in the solution. Recently \cite{price2017renormalized}, we have developed a novel expansion of the memory term (dubbed the ``complete memory approximation") which also assumes long memory but avoids all the simplifying approximations applied before. Thus, it can incorporate all the complex effects present in the Euler dynamics. 

In the current work we present results for renormalized MZ reduced models of the 3D Euler equations stemming from the complete memory approximation. To the best of our knowledge this is the first time-dependent perturbative renormalization approach for 3D Euler which includes all the complex effects present in the Euler dynamics. For the Taylor-Green initial condition, the renormalized coefficients in the reduced order models decay algebraically with time and resolution. The renormalized reduced models are stable and we evolve them for long times. We use the predicted evolution to estimate several quantities related to the occurrence of a finite-time singularity (blow-up).  We find that our results are indeed consistent with a finite-time singularity. 

This paper is organized in the following manner. In Section 2, we present a brief overview of the MZ formalism and constructing reduced order models from it by employing the complete memory approximation. In Section 3, we demonstrate the application of this approach to Euler's equations including the process of identifying effective renormalization coefficients. In Section 4, we present the results of simulations of these reduced order models. In Section 5, we discuss these results and comment upon the question of finite-time singularities in Euler's equations.

\section{The Mori-Zwanzig Formalism}

Previous work \cite{price2017renormalized} includes a comprehensive overview of the Mori-Zwanzig formalism and construction of reduced order models from it by way of the complete memory approximation. Here we present an abridged version. Consider a system of autonomous ordinary differential equations
\begin{equation}
\frac{d \mathbf{u}(t)}{dt} = \mathbf{R}(\mathbf{u})\label{orig}
\end{equation}
augmented with an initial condition $\mathbf{u}(0)=\mathbf{u}^0$. Let $\mathbf{u}(t) = \{u_k(t)\}$, $k\in F\cup G$. We separate $\mathbf{u}(t)$ into resolved variables $\hat{\mathbf{u}}=\{u_i(t)\}$, $i\in F$ and unresolved variables $\tilde{\mathbf{u}}=\{u_j(t)\}$, $j\in G$ where $F$ and $G$ are disjoint. Let $R_k(\mathbf{u})$ be the $k$th entry in the vector-valued function $\mathbf{R}(\mathbf{u})$. We can transform this nonlinear system of ODEs (\ref{orig}) into a linear system of PDEs by way of the Liouvillian operator \cite{chorin2000optimal,chorin2002optimal}, also known as the generator of the Koopman operator \cite{koopman1931hamiltonian}:
\begin{equation}
\mathcal{L}=\sum_{k\in F\cup G} R_k(\mathbf{u}^0)\frac{\partial}{\partial u_k^0}.\label{L}
\end{equation}
It can be shown that $\phi_k(\mathbf{u}^0,t) = u_k(\mathbf{u}^0,t)$ satisfies the PDE:
\begin{equation}
\frac{\partial \phi_k(\mathbf{u}^0,t)}{\partial t} = \mathcal{L}\phi_k,\;\; \phi_k(\mathbf{u}^0,0)=u_k^0.\label{liouville}
\end{equation} In essence, we are constructing the partial differential equation for which our original ODE system defines the backward characteristic curves. Equation (\ref{liouville}), however, is linear. Using semigroup notation, which states $\phi_k(\mathbf{u}^0,t) = e^{t\mathcal{L}}u_k^0$, we can write
\begin{equation}
\frac{du_k(t)}{dt} = \frac{\partial \phi_k(\mathbf{u}^0,t)}{\partial t} = \frac{\partial}{\partial t}e^{t\mathcal{L}}u_k^0 = e^{t\mathcal{L}}\mathcal{L}u_k^0.\label{semigroup}
\end{equation} 
Consider the space of functions that depend upon $\mathbf{u}^0$. Let $P$ be an orthogonal projection onto the subspace of functions depending only on the resolved variables $\hat{\mathbf{u}}^0$. For example, $Pf$ might be the conditional expectation of $f$ given the resolved variables and an assumed joint density. Let $Q=I-P$. Then, we can decompose the evolution operator $e^{t\mathcal{L}}$ using Dyson's formula,
\begin{equation}
e^{t\mathcal{L}}=e^{tQ\mathcal{L}}+\int_0^t e^{(t-s)\mathcal{L}}P\mathcal{L}e^{sQ\mathcal{L}}\,\mathrm{d}s\label{dyson}
\end{equation}into:
\begin{equation}
\frac{du_k}{dt} = e^{t\mathcal{L}}P\mathcal{L}u_k^0 +e^{tQ\mathcal{L}}Q\mathcal{L}u_k^0 + \int_0^t e^{(t-s)\mathcal{L}}P\mathcal{L}e^{sQ\mathcal{L}}Q\mathcal{L}u_k^0\,\mathrm{d}s.\label{MZ}
\end{equation}This is the Mori-Zwanzig identity. It is simply a rewritten version of the original dynamics. The first term on the right hand side in (\ref{MZ}) is called the Markov term, because it depends only on the instantaneous values of the resolved variables. The second term is called `noise' and the third is called `memory'. We again project the dynamics:
\begin{equation}
\frac{dPu_k}{dt} = Pe^{t\mathcal{L}}P\mathcal{L}u_k^0+P\int_0^te^{(t-s)\mathcal{L}}P\mathcal{L}e^{sQ\mathcal{L}}Q\mathcal{L}u_k^0\,\mathrm{d}s.\label{reduced_MZ}
\end{equation} For $k\in F$, (\ref{reduced_MZ}) describes the projected dynamics of the resolved variables. The system is not closed, however, due to the presence of the orthogonal dynamics operator $e^{sQ\mathcal{L}}$ in the memory term. In order to simulate the dynamics of (\ref{reduced_MZ}) exactly, one needs to evaluate the second term which requires the dynamics of the unresolved variables. 

One key fact must be understood: reducing a large system to one of comparatively fewer variables necessarily introduces a memory term encoding the interplay between the unresolved and resolved variables. Dropping both the noise and memory terms and simulating only the ``average'' dynamics (the Markov term) may not accurately reflect the dynamics of the resolved variables in the full simulation. Any multiscale dynamical model must in some way approximate or compute the memory term, or argue convincingly why the memory term is negligible. In a previous work, it was shown that even when the memory term is small in magnitude, neglecting it leads to inaccurate simulations \cite{price2017renormalized}.

\subsection{Renormalized Memory Approximations}

Define the Markov term as
\begin{equation}
Pe^{t\mathcal{L}}P\mathcal{L}u_k^0 = R_k^0(\hat{\mathbf{u}}),
\end{equation}
and the memory term as
\begin{equation}
\mathcal{M}_k = P\int_0^te^{(t-s)\mathcal{L}}P\mathcal{L}e^{sQ\mathcal{L}}Q\mathcal{L}u_k^0\,\mathrm{d}s.
\end{equation}
The simplest possible approximation of the memory integral is to assume the integrand is constant. In this case, the memory integral becomes:
\begin{equation}
\mathcal{M}_k \approx tPe^{t\mathcal{L}}P\mathcal{L}Q\mathcal{L}u_k^0.
\end{equation}We apply $P\mathcal{L}Q\mathcal{L}$ to the initial condition. The result is an expression that depends upon only the initial condition of the resolved modes (because the final operator applied before the evolution operator is a projection operator). Therefore, this term can be expressed in terms of only the trajectories of the resolved variables, as desired. This model is called the $t$-model and it has been used to successfully construct reduced order models for a variety of problems \cite{chorin2002optimal,stinis2012numerical,chorin2007problem,hald2007optimal,bernstein2007optimal,chandy2009t}. 

\subsubsection{Complete Memory Approximation}

We would like to improve upon the accuracy of the $t$-model. It is our hope to construct a series representation of $\mathcal{M}_k$ in powers of $t$. We begin by rewriting the memory term using the definitions of $e^{-s\mathcal{L}}$ and $e^{sQ\mathcal{L}}$ and then computing the integral termwise:
\begin{align}
\mathcal{M}_k
=&Pe^{t\mathcal{L}}\int_0^t \left(\sum_{i=0}^{\infty}\frac{(-1)^is^i}{i!}\mathcal{L}^i\right)P\mathcal{L}\left(\sum_{j=0}^\infty \frac{s^j}{j!}(Q\mathcal{L})^j\right)Q\mathcal{L}u_k^0\,\mathrm{d}s\notag\\
=&Pe^{t\mathcal{L}}\left(\sum_{i=0}^\infty \sum_{j=0}^\infty \frac{(-1)^it^{i+j+1}}{i!j!(i+j+1)}\mathcal{L}^iP\mathcal{L}(Q\mathcal{L})^{j}Q\mathcal{L}u_k^0\right).\label{eq:novel}
\end{align}We assume the integrand is sufficiently smooth that we can interchange the order of the integral and the infinite sums. Now, consider the full expansion of the memory kernel \eqref{eq:novel}. We begin to analyze this formulation of the memory by writing the first few terms arranged by powers of $t$:
\begin{align}
\mathcal{M}_k =&tPe^{t\mathcal{L}}\left[P\mathcal{L}Q\mathcal{L}\right]u_k^0 - \frac{t^2}{2}Pe^{t\mathcal{L}}\left[\mathcal{L}P\mathcal{L}Q\mathcal{L}-P\mathcal{L}Q\mathcal{L}Q\mathcal{L}\right]u_k^0+O(t^3) \label{unclosed}
\end{align}The $O(t)$ term is the $t$-model once again. The $O(t^2)$ term, presents a new problem. The second term in it can be computed in a manner similar to the $t$-model, but the first term is not projected prior to its evolution. $\mathcal{L}f(\mathbf{u}^0)$ is a function of all modes, not just the unresolved ones. If we wish to evolve a term of this form forward in time, we would need to evolve forward in time a quantity that depends upon unresolved modes, necessitating knowledge of the dynamics of the unresolved modes. This makes it impossible to compute as part of a reduced order model except in very special cases. 

To close the model in the resolved variables we construct an additional reduced order model for the problem term. This term is $Pe^{t\mathcal{L}}\mathcal{L}P\mathcal{L}Q\mathcal{L}u_k^0$. First note that
\begin{equation}
Pe^{t\mathcal{L}}\mathcal{L}P\mathcal{L}Q\mathcal{L}u_k^0=\frac{\partial}{\partial t} Pe^{t\mathcal{L}}P\mathcal{L}Q\mathcal{L}u_k^0.
\end{equation}That is, it is the derivative of the $t$-model term itself. Now consider a reduced order model for this derivative under the Mori-Zwanzig formalism again:
\begin{align}
\frac{\partial}{\partial t} &Pe^{t\mathcal{L}}P\mathcal{L}Q\mathcal{L}u_k^0 =Pe^{t\mathcal{L}}P\mathcal{L}P\mathcal{L}Q\mathcal{L}u_k^0+P\int_0^t e^{(t-s)\mathcal{L}}P\mathcal{L}e^{sQ\mathcal{L}}Q\mathcal{L}P\mathcal{L}Q\mathcal{L}u_k^0\,\mathrm{d}s\notag\\
=&Pe^{t\mathcal{L}}P\mathcal{L}P\mathcal{L}Q\mathcal{L}u_k^0 + Pe^{t\mathcal{L}}\left(\sum_{i=0}^\infty\sum_{j=0}^\infty \frac{(-1)^it^{i+j+1}}{i!j!(i+j+1)}\mathcal{L}^iP\mathcal{L}(Q\mathcal{L})^{j+1}P\mathcal{L}Q\mathcal{L}u_k^0\right).\label{tmodeloftmodel}
\end{align}
If we replace $Pe^{t\mathcal{L}}\mathcal{L}P\mathcal{L}Q\mathcal{L}u_k^0$ in \eqref{unclosed} with \eqref{tmodeloftmodel}, only the first term of  \eqref{tmodeloftmodel} contributes at the $O(t^2)$ level. We find 
\begin{align}
\mathcal{M}_k =&tPe^{t\mathcal{L}}\left[P\mathcal{L}Q\mathcal{L}\right]u_k^0 - \frac{t^2}{2}Pe^{t\mathcal{L}}P\mathcal{L}\left[P\mathcal{L}-Q\mathcal{L}\right]Q\mathcal{L}u_k^0+O(t^3) \label{closed}
\end{align}
where all the terms are now projected prior to evolution and so involve {\it only} resolved variables. 

Many of the higher order terms in \eqref{eq:novel} contain a leading $\mathcal{L}$ like the ``problem term'' at $O(t^2)$ discussed above. In each case, we can construct a reduced order model for the problem term and approximate it by expanding the memory term in a series. The terms in this series will also include leading $\mathcal{L}$ terms, but we can simply repeat our procedure indefinitely. In this manner, we can construct an approximation for \eqref{eq:novel} in which every term has a leading $P$ before the evolution operator is applied. The resulting series is written as:
\begin{equation}
\mathcal{M}_k = \sum_{i=1}^\infty \frac{(-1)^{i+1}t^i}{i!} R_k^i(\hat{\mathbf{u}}).
\end{equation}Different approximation schemes can be constructed by truncating this series at different terms. As discussed above, the $O(t)$ term corresponds to the $t$-model\begin{equation}
R_k^1(\hat{\mathbf{u}}) = Pe^{t\mathcal{L}}P\mathcal{L}Q\mathcal{L}u_k^0,\label{tmodel}
\end{equation}
and the $O(t^2)$ term is:
\begin{equation}
R_k^2(\hat{\mathbf{u}}) = Pe^{t\mathcal{L}}P\mathcal{L}\left[P\mathcal{L}-Q\mathcal{L}\right]Q\mathcal{L}u_k^0.\label{t2}
\end{equation}By grouping terms in the series \eqref{eq:novel} in powers of $t$ and using the technique described above to further expand problem terms, we can uniquely define $R_k^i(\hat{\mathbf{u}})$ for any positive integer $i$. We automated this process in a Mathematica notebook, which is available in the \texttt{Renormalized\_Mori\_Zwanzig} git repository \cite{MZrepos}. The $O(t^3)$ term is:
\begin{equation}
R_k^3(\hat{\mathbf{u}}) = Pe^{t\mathcal{L}}P\mathcal{L}[P\mathcal{L}P\mathcal{L}-2P\mathcal{L}Q\mathcal{L}-2Q\mathcal{L}P\mathcal{L}+Q\mathcal{L}Q\mathcal{L}]Q\mathcal{L}u_k^0.\label{t3}
\end{equation}
Finally, the $O(t^4)$ term is:
\begin{align}
R_k^4(\hat{\mathbf{u}}) = Pe^{t\mathcal{L}}P\mathcal{L}\bigg[&P\mathcal{L}P\mathcal{L}P\mathcal{L}-3P\mathcal{L}P\mathcal{L}Q\mathcal{L}-5P\mathcal{L}Q\mathcal{L}P\mathcal{L}-3Q\mathcal{L}P\mathcal{L}P\mathcal{L}\notag\\
&+3P\mathcal{L}Q\mathcal{L}Q\mathcal{L}+5Q\mathcal{L}P\mathcal{L}Q\mathcal{L}+3Q\mathcal{L}Q\mathcal{L}P\mathcal{L}-Q\mathcal{L}Q\mathcal{L}Q\mathcal{L}\bigg]Q\mathcal{L}u_k^0.\label{t4}
\end{align}

\subsubsection{Renormalization}

The complete memory approximation framework provides a series representation of the memory integral. Different ROMs can be created by truncating the series at different terms. In this case, the differential equation for a resolved mode is:
\begin{equation}
\frac{dPu_k}{dt} = R_k^0(\hat{\mathbf{u}}) + \sum_{i=1}^n \frac{(-1)^{i+1}t^i}{i!} R_k^i(\hat{\mathbf{u}}).,\label{nonrenormalized}
\end{equation}When applied to Euler's equations, the resulting ROMs are unstable. 

For a related memory approximation method, it was found that \emph{renormalization} rendered the reduced order models for Euler's equations stable \cite{stinis2013renormalized,stinis2015renormalized}. We attach additional coefficients to each term in the series, such that the terms represent an \emph{effective memory}, given knowledge only of the resolved modes \cite{georgi1993}. The evolution equation for a reduced variable then becomes
\begin{equation}
\frac{dPu_k}{dt} = R_k^0(\hat{\mathbf{u}}) + \sum_{i=1}^n \alpha_i(t)t^i R_k^i(\hat{\mathbf{u}}).\label{renormalized}
\end{equation}Here, we allow the renormalization coefficients $\alpha_i(t)$ to be time dependent. This gives us the flexibility to allow the functional form of the effective memory to be dynamic if necessary. These coefficients must be chosen in a way that captures information we know about the memory term. We will detail how we selected coefficients in Section \ref{renorm}.

\section{Reduced Order Models of Euler's Equations in Three Dimensions}\label{euler}

The three dimensional Euler's equations are given by:
\begin{equation}
\mathbf{u}_t + \mathbf{u}\cdot\nabla\mathbf{u} = -\nabla p,\quad \nabla\cdot \mathbf{u} = 0\label{Euler}
\end{equation}where $\mathbf{u}(x,t)$ is the three-dimensional velocity field and $p$ is the pressure. We restrict attention to periodic boundary conditions. Also, let the initial condition be $\mathbf{u}(x,0) = \mathbf{v}^0$. We will use the Taylor-Green initial condition:
\begin{equation}
\mathbf{v}^0 = \begin{bmatrix}
\sin(x)\cos(y)\cos(z)\\
-\cos(x)\sin(y)\cos(z)\\
0
\end{bmatrix}.\label{taylorgreen}
\end{equation}The Taylor-Green initial condition is very smooth. We are interested in studying the cascade of energy into higher frequency modes as time evolves. Because $\mathbf{u}$ is periodic in all three dimensions, we will write it as Fourier series:
\begin{equation}
\mathbf{u}(x,t) = \sum_{\mathbf{k}} \mathbf{u}_{\mathbf{k}}(t) e^{i\mathbf{k}\cdot\mathbf{x}}\label{eq:fourier}
\end{equation}where $\mathbf{k}$ is a three-dimensional wavevector, and the sum is over all possible integer-valued wavevectors. The evolution of a Fourier mode is:
\begin{equation}
\frac{d\mathbf{u}_{\mathbf{k}}}{dt} =\mathbf{R}_\mathbf{k}(\mathbf{u}) =  -i\sum_{\mathbf{p}+\mathbf{q}=\mathbf{k}}\mathbf{k}\cdot\mathbf{u}_{\mathbf{p}}A_{\mathbf{k}}\mathbf{u}_{\mathbf{q}}\label{evolution}
\end{equation}where
\begin{equation}
A_k = I - \frac{\mathbf{k}\mathbf{k}^T}{|\mathbf{k}|^2}
\end{equation}
is the incompressibility projection operator \cite{doering1995applied}. We will use Matlab's built-in integrator \texttt{ode45} to solve the described differential equations. This uses a version of Runge-Kutta-Fehlberg with adaptive stepsize selection. We set the initial step as $10^{-3}$ because the very small initial rate of change causes the algorithm to choose an overly ambitious starting step. We also set the maximum relative error to $10^{-10}$.

\subsection{Reduced models}

Let $\mathbf{u}(x,t)$ be the solution to Euler's equations. Consider the Fourier components $\mathbf{u}_{\mathbf{k}}(t)$, where $\mathbf{k}\in F\cup G$. Let $F$ be the set of resolved modes. That is $F = \{\mathbf{k} \in [-N,N-1]^3\}$. Let  $F\cup G = \{\mathbf{k} \in [-M,M-1]^3\}$. Define $\hat{\mathbf{u}} = \{\mathbf{u}_{\mathbf{k}}\;\; |\;\; \mathbf{k}\in F\}$ and $\tilde{\mathbf{u}} = \{\mathbf{u}_{\mathbf{k}}\;\; | \;\;\mathbf{k} \in G\}$.

Following the Mori-Zwanzig formalism, we define
\[\mathcal{L} = \sum_{\mathbf{k}\in F\cup G} \begin{bmatrix}R^x_{\mathbf{k}}(\mathbf{u}^0)\frac{\partial}{\partial u_{\mathbf{k}}^{0,x}} & 0 & 0\\
0 &R^y_{\mathbf{k}}(\mathbf{u}^0)\frac{\partial}{\partial u_{\mathbf{k}}^{0,y}} & 0\\
0 & 0 &R^z_{\mathbf{k}}(\mathbf{u}^0)\frac{\partial}{\partial u_{\mathbf{k}}^{0,z}} \end{bmatrix}
\]where $R_{\mathbf{k}}^i$ is the component of $\mathbf{R}_{\mathbf{k}}$ in the $i$th direction and $u_{\mathbf{k}}^{0,i}$ is the component of the initial condition of $\mathbf{u}_{\mathbf{k}}$ in the $i$th direction. Observe that:
\[\mathcal{L}\mathbf{u}_{\mathbf{k}}^0 = \mathbf{R}_{\mathbf{k}}(\mathbf{u}^0)
\]as desired. Therefore, we can see
\[\mathcal{L}\mathbf{u}_{\mathbf{k}}^0 = -i\sum_{\substack{\mathbf{p}+\mathbf{q}=\mathbf{k}\\ \mathbf{p},\mathbf{q},\mathbf{k}\in F\cup G}}\mathbf{k}\cdot \mathbf{u}^0_{\mathbf{p}}A_{\mathbf{k}}\mathbf{u}^0_{\mathbf{q}}.
\]

We define the convolution $\mathbf{C}$ of two vectors. Let $\mathbf{C}_{\mathbf{k}}$ be the component of $\mathbf{C}$ corresponding to the wavevector $\mathbf{k}$:
\begin{equation}\mathbf{C}_{\mathbf{k}}(\mathbf{v},\mathbf{w}) = -i\sum_{\substack{\mathbf{p}+\mathbf{q}=\mathbf{k}\\\mathbf{p},\mathbf{q} \in F\cup G}}\mathbf{k}\cdot\mathbf{v}_{\mathbf{p}} A_{\mathbf{k}} \mathbf{w}_{\mathbf{q}}.
\end{equation}
Under this definition, we can see that $\mathcal{L}\mathbf{u}_{\mathbf{k}}^0= \mathbf{C}_{\mathbf{k}}(\mathbf{u}^0,\mathbf{u}^0)$, or under the Mori-Zwanzig formalism that \begin{equation}\frac{d\mathbf{u}_{\mathbf{k}}}{dt} = e^{t\mathcal{L}}\mathbf{C}_{\mathbf{k}}(\mathbf{u}^0,\mathbf{u}^0) = \mathbf{C}_{\mathbf{k}}(\mathbf{u},\mathbf{u}).\label{full_Euler}
\end{equation}For this to be the exact full solution, we would need $M = \infty$. For our implementation, $M$ will necessarily be the maximal mode we choose to retain.

We now define the projector $P$ as:
\[P f(\mathbf{u}) = P f(\hat{\mathbf{u}}^0,\tilde{\mathbf{u}}^0) = f(\hat{\mathbf{u}}^0,0).
\]That is, it simply sets all unresolved modes to zero. Let $Q = I - P$. Then we can begin constructing the terms of a complete memory approximation:
\begin{equation}\frac{dP\mathbf{u}_{\mathbf{k}}}{dt} = \mathbf{R}_{\mathbf{k}}^{0}(\hat{\mathbf{u}})+\sum_{i=1}^4 \alpha_i(t)t^i\mathbf{R}_{\mathbf{k}}^{i}(\hat{\mathbf{u}}).
\end{equation}

\subsubsection{Markov Term}

The Markov term is
\[\mathbf{R}^0_{\mathbf{k}}(\hat{\mathbf{u}}) = Pe^{t\mathcal{L}}P\mathcal{L}\mathbf{u}_{\mathbf{k}}^0.
\]Using the definitions above, we find:
\begin{align*}\mathbf{R}^0_{\mathbf{k}}(\hat{\mathbf{u}}) =& Pe^{t\mathcal{L}}P\mathbf{C}_{\mathbf{k}}(\mathbf{u}^0,\mathbf{u}^0)
= Pe^{t\mathcal{L}}\mathbf{C}_{\mathbf{k}}(\hat{\mathbf{u}}^0,\hat{\mathbf{u}}^0)
= \mathbf{C}_{\mathbf{k}}(\hat{\mathbf{u}},\hat{\mathbf{u}}).
\end{align*}Here, we used the fact that $P$ sets the unresolved modes to zero. Through a slight abuse of notation, consider $\hat{\mathbf{u}}$ to be the array of $\mathbf{u}$ where all $\mathbf{u}_{\mathbf{j}} = 0$ for $\mathbf{j}\in G$ (and similarly for $\tilde{\mathbf{u}}$).

\subsubsection{First-order term}

The $t$-model is
\[\mathbf{R}^1_{\mathbf{k}}(\hat{\mathbf{u}}) = Pe^{t\mathcal{L}}P\mathcal{L}Q\mathcal{L}\mathbf{u}_{\mathbf{k}}^0.
\]We can simplify this as:
\begin{align*}
\mathbf{R}^1_{\mathbf{k}}(\hat{\mathbf{u}}) =& Pe^{t\mathcal{L}}P\mathcal{L}Q\mathcal{L}\mathbf{u}_{\mathbf{k}}^0
= Pe^{t\mathcal{L}}P\mathcal{L}[\mathcal{L}-P\mathcal{L}]\mathbf{u}_{\mathbf{k}}^0\\
=& Pe^{t\mathcal{L}}P\mathcal{L} [\mathbf{C}_{\mathbf{k}}(\hat{\mathbf{u}}^0,\tilde{\mathbf{u}}^0)+\mathbf{C}_{\mathbf{k}}(\tilde{\mathbf{u}}^0,\hat{\mathbf{u}}^0)+\mathbf{C}_{\mathbf{k}}(\tilde{\mathbf{u}}^0,\tilde{\mathbf{u}}^0)]
\end{align*}When applied to a convolution sum, $\mathcal{L}$ obeys the product rule \[\mathcal{L}\mathbf{C}_\mathbf{k}(\mathbf{v},\mathbf{w}) = \mathbf{C}_\mathbf{k}(\mathcal{L}\mathbf{v},\mathbf{w})+\mathbf{C}_\mathbf{k}(\mathbf{v},\mathcal{L}\mathbf{w}).\]Define the following arrays. Let $\hat{\mathbf{C}}(\mathbf{v},\mathbf{w})$ be the convolution of $\mathbf{v}$ and $\mathbf{w}$, but set all indices corresponding to $\mathbf{j}\in G$ to zero.  Let $\tilde{\mathbf{C}}(\mathbf{v},\mathbf{w})$ be the convolution of $\mathbf{v}$ and $\mathbf{w}$, but set all indices corresponding to $\mathbf{i}\in F$ to zero. Thus, $\mathbf{C}(\mathbf{v},\mathbf{w}) = \hat{\mathbf{C}}(\mathbf{v},\mathbf{w})+\tilde{\mathbf{C}}(\mathbf{v},\mathbf{w})$. Observe that, under this definition, $\mathcal{L}\hat{\mathbf{u}}^0 = \hat{\mathbf{C}}(\mathbf{u}^0,\mathbf{u}^0)$ and $\mathcal{L}\tilde{\mathbf{u}}^0 = \tilde{\mathbf{C}}(\mathbf{u}^0,\mathbf{u}^0)$. Thus,
\begin{align*}
\mathbf{R}^1_{\mathbf{k}}(\hat{\mathbf{u}})
=& Pe^{t\mathcal{L}}P[\mathbf{C}_{\mathbf{k}}(\hat{\mathbf{C}}(\mathbf{u}^0,\mathbf{u}^0),\tilde{\mathbf{u}}^0)+\mathbf{C}_{\mathbf{k}}(\hat{\mathbf{u}}^0,\tilde{\mathbf{C}}(\mathbf{u}^0,\mathbf{u}^0))\\
&\qquad\quad \;+\mathbf{C}_{\mathbf{k}}(\tilde{\mathbf{C}}(\mathbf{u}^0,\mathbf{u}^0),\hat{\mathbf{u}}^0)+\mathbf{C}_{\mathbf{k}}(\tilde{\mathbf{u}}^0,\hat{\mathbf{C}}(\mathbf{u}^0,\mathbf{u}^0))\\
&\qquad\quad \;+\mathbf{C}_{\mathbf{k}}(\tilde{\mathbf{C}}(\mathbf{u}^0,\mathbf{u}^0),\tilde{\mathbf{u}}^0)+\mathbf{C}_{\mathbf{k}}(\tilde{\mathbf{u}}^0,\tilde{\mathbf{C}}(\mathbf{u}^0,\mathbf{u}^0))].
\end{align*}In order to apply the projector, observe that $P\mathbf{C}_{\mathbf{k}}(\mathbf{v},\mathbf{w}) = \mathbf{C}_{\mathbf{k}}(P\mathbf{v},P\mathbf{w})$. Note that $P\mathbf{u}^0 = \hat{\mathbf{u}}^0$, $P\hat{\mathbf{u}}^0 = \hat{\mathbf{u}}^0$, and $P\tilde{\mathbf{u}}^0 = 0$. Finally, note that $\mathbf{C}_{\mathbf{k}}(\mathbf{v},0) = \mathbf{C}_{\mathbf{k}}(0,\mathbf{w}) = 0$. This yields:
\begin{align*}
\mathbf{R}^1_{\mathbf{k}}(\hat{\mathbf{u}}) =& Pe^{t\mathcal{L}}P[\mathbf{C}_{\mathbf{k}}(\hat{\mathbf{C}}(\mathbf{u}^0,\mathbf{u}^0),\tilde{\mathbf{u}}^0)+\mathbf{C}_{\mathbf{k}}(\hat{\mathbf{u}}^0,\tilde{\mathbf{C}}(\mathbf{u}^0,\mathbf{u}^0))\\
&\qquad\quad \;+\mathbf{C}_{\mathbf{k}}(\tilde{\mathbf{C}}(\mathbf{u}^0,\mathbf{u}^0),\hat{\mathbf{u}}^0)+\mathbf{C}_{\mathbf{k}}(\tilde{\mathbf{u}}^0,\hat{\mathbf{C}}(\mathbf{u}^0,\mathbf{u}^0))\\
&\qquad\quad \;+\mathbf{C}_{\mathbf{k}}(\tilde{\mathbf{C}}(\mathbf{u}^0,\mathbf{u}^0),\tilde{\mathbf{u}}^0)+\mathbf{C}_{\mathbf{k}}(\tilde{\mathbf{u}}^0,\tilde{\mathbf{C}}(\mathbf{u}^0,\mathbf{u}^0))]\\
=& Pe^{t\mathcal{L}}[\mathbf{C}_{\mathbf{k}}(\hat{\mathbf{u}}^0,\tilde{\mathbf{C}}(\hat{\mathbf{u}}^0,\hat{\mathbf{u}}^0)) + \mathbf{C}_{\mathbf{k}}(\tilde{\mathbf{C}}(\hat{\mathbf{u}}^0,\hat{\mathbf{u}}^0),\hat{\mathbf{u}}^0)]\\
=&\mathbf{C}_{\mathbf{k}}(\hat{\mathbf{u}},\tilde{\mathbf{C}}(\hat{\mathbf{u}},\hat{\mathbf{u}})) + \mathbf{C}_{\mathbf{k}}(\tilde{\mathbf{C}}(\hat{\mathbf{u}},\hat{\mathbf{u}}),\hat{\mathbf{u}}).
\end{align*}This term involves convolutions of the same two terms in both permutations. It is useful to define a function:
\begin{align*}
\mathbf{D}(\mathbf{v},\mathbf{w}) = \mathbf{C}(\mathbf{v},\mathbf{w}) + \mathbf{C}(\mathbf{w},\mathbf{v})
\end{align*}and the related functions $\mathbf{D}_{\mathbf{k}}$, $\hat{\mathbf{D}}$, and $\tilde{\mathbf{D}}$ defined in terms of the equivalent convolutions. With this notation, we can write the $t$-model term in a single expression:
\[\mathbf{R}^1_{\mathbf{k}}(\hat{\mathbf{u}}) = \mathbf{D}_{\mathbf{k}}(\hat{\mathbf{u}},\tilde{\mathbf{C}}(\hat{\mathbf{u}},\hat{\mathbf{u}}))
\]

\subsubsection{Second-order term}

The set of rules derived in the previous section will allow us to proceed to higher terms. The second order term is:
\begin{align*}\mathbf{R}^2_{\mathbf{k}}(\hat{\mathbf{u}}) =& Pe^{t\mathcal{L}}P\mathcal{L}[P\mathcal{L}-Q\mathcal{L}]Q\mathcal{L}\mathbf{u}_{\mathbf{k}}^0.
\end{align*}First note that, with our newly defined $\mathbf{D}$ function, \[Q\mathcal{L}\mathbf{u}_{\mathbf{k}}^0 = \mathbf{D}_{\mathbf{k}}(\hat{\mathbf{u}},\tilde{\mathbf{u}}) + \mathbf{C}_{\mathbf{k}}(\tilde{\mathbf{u}},\tilde{\mathbf{u}}).\] $\mathcal{L}$ operates upon $\mathbf{D}$ in the same manner it did upon $\mathbf{C}$. The projector $P$ when applied to $\mathbf{D}$ similarly is applied to each term within the expression. Starting from this, we derive an expression for the $t^2$-term:
\begin{align*}
\mathbf{R}^2_{\mathbf{k}}(\hat{\mathbf{u}}) =& Pe^{t\mathcal{L}}P\mathcal{L}(2P\mathcal{L}-\mathcal{L})[\mathbf{D}_{\mathbf{k}}(\hat{\mathbf{u}}^0,\tilde{\mathbf{u}}^0)+\mathbf{C}_{\mathbf{k}}(\tilde{\mathbf{u}}^0,\tilde{\mathbf{u}}^0)]\\
=&Pe^{t\mathcal{L}}P\mathcal{L}[2\mathbf{D}_{\mathbf{k}}(\hat{\mathbf{u}}^0,\tilde{\mathbf{C}}(\hat{\mathbf{u}}^0,\hat{\mathbf{u}}^0))-\mathbf{D}_{\mathbf{k}}(\hat{\mathbf{C}}(\mathbf{u}^0,\mathbf{u}^0),\tilde{\mathbf{u}}^0)\\
&\qquad\qquad\quad -\mathbf{D}_{\mathbf{k}}(\hat{\mathbf{u}}^0,\tilde{\mathbf{C}}(\mathbf{u}^0,\mathbf{u}^0))-\mathbf{D}_{\mathbf{k}}(\tilde{\mathbf{u}}^0,\tilde{\mathbf{C}}(\mathbf{u}^0,\mathbf{u}^0))]\\
=&Pe^{t\mathcal{L}}[\mathbf{D}_{\mathbf{k}}(\hat{\mathbf{u}}^0,\tilde{\mathbf{D}}(\hat{\mathbf{C}}(\hat{\mathbf{u}}^0,\hat{\mathbf{u}}^0)-\tilde{\mathbf{C}}(\hat{\mathbf{u}}^0,\hat{\mathbf{u}}^0),\hat{\mathbf{u}}^0))-\mathbf{D}_{\mathbf{k}}(\tilde{\mathbf{C}}(\hat{\mathbf{u}}^0,\hat{\mathbf{u}}^0),\tilde{\mathbf{C}}(\hat{\mathbf{u}}^0,\hat{\mathbf{u}}^0))]\\
=&\mathbf{D}_{\mathbf{k}}(\hat{\mathbf{u}},\tilde{\mathbf{D}}(\hat{\mathbf{C}}(\hat{\mathbf{u}},\hat{\mathbf{u}})-\tilde{\mathbf{C}}(\hat{\mathbf{u}},\hat{\mathbf{u}}),\hat{\mathbf{u}}))-\mathbf{D}_{\mathbf{k}}(\tilde{\mathbf{C}}(\hat{\mathbf{u}},\hat{\mathbf{u}}),\tilde{\mathbf{C}}(\hat{\mathbf{u}},\hat{\mathbf{u}})).
\end{align*}

\subsubsection{Higher ordered terms}

The rules described in the previous sections can be automated in a symbolic notebook. We used a notebook (included in the git repository \cite{MZrepos}) to compute third and fourth order terms in the complete memory approximation for Euler's equations. These terms and details of this derivation are described in Appendix \ref{sec:higherorder}. The terms are quite complicated and are unlikely to be proposed by a mathematical modeler. However, we find that these terms lead to stable reduced order models with significant structure in the renormalization coefficients. For this reason, we propose that the complete memory approximation represents a ``natural'' reduced order model derived from the full equations themselves.

\subsection{Renormalization coefficients}\label{renorm}

With these terms computed, we can express a renormalized reduced order model as:
\[\frac{dP\mathbf{u}_{\mathbf{k}}}{dt} =\mathbf{R}_{\mathbf{k}}^{0}(\hat{\mathbf{u}}) +  \sum_{i=1}^n \alpha_i(t)t^i\mathbf{R}_{\mathbf{k}}^{i}(\hat{\mathbf{u}})\]for $n = 1,2,3,4$. If we do not renormalize and $\alpha_i(t) = \frac{(-1)^{i+1}}{i!}$, we find that the simulations are unstable for all except the $t$-model alone. Instead, we will choose renormalization coefficients to stabilize the models. In past work, it was found that ROMs of Burgers' equation and Euler's equations were stabilized with constant renormalization coefficients \cite{stinis2013renormalized,stinis2015renormalized} and ROMs of the Korteweg-de Vries equation were stabilized with algebraically decaying renormalization coefficients \cite{price2017renormalized}. We will consider both cases as possible ansatzes. To be explicit,
\begin{align}
\alpha_i(t) &= a_it^{-i}\\
\alpha_i'(t) &= a_i'.
\end{align}
We will use the rates of change of the energy in each resolved mode as the quantities we attempt to match in the renormalization process. This choice is reasonable because it is known that energy moves from low-frequency modes to high-frequency modes as Euler's equations are evolved, and the Markov term is incapable of capturing this, since it conserves energy in the resolved modes. Thus, it makes an excellent heuristic for the effectiveness of a memory approximation. The energy of a mode is defined as:
\begin{equation}
E_{\mathbf{k}}(t) = |\mathbf{u}_{\mathbf{k}}|^2.
\end{equation}The rate of change of the energy in a particular mode in the full model is:
\begin{equation}
\Delta E_{\mathbf{k}}(t) = \mathbf{R}_{\mathbf{k}}(\mathbf{u})\cdot\overline{\mathbf{u}}_{\mathbf{k}}+\mathbf{u}_{\mathbf{k}}\cdot\overline{\mathbf{R}}_{\mathbf{k}}(\mathbf{u}).
\end{equation}In a reduced order model, each term in the series has its own contribution to the energy derivative:
\begin{equation}
\Delta E_{\mathbf{k}}^i(t) = \mathbf{R}^i_{\mathbf{k}}(\hat{\mathbf{u}})\cdot\overline{\mathbf{u}}_{\mathbf{k}}+\mathbf{u}_{\mathbf{k}}\cdot\overline{\mathbf{R}}^i_{\mathbf{k}}(\hat{\mathbf{u}}).
\end{equation}Given an exact energy derivative, our renormalization coefficients will be chosen to minimize the difference between $\Delta E_{\mathbf{k}}$ and 
\[\Delta \hat{E}_{\mathbf{k}}^n(t) = \Delta E_{\mathbf{k}}^0(t)+\sum_{i=1}^n \alpha_i(t)t^i\Delta E^i_{\mathbf{k}}(t).
\]

We will renormalize against data $\Delta E_{\mathbf{k}}$ produced by a full model that we trust has not yet become unresolved. We use the Markov model to simulate this ``full'' system of size $F \cup G = \{\mathbf{k} \;\;|\;\; \mathbf{k}\in [-M,M-1]^3\}$ up to time $T$. This produces a time series $\mathbf{u}(t)$ for $t = 0,\dots,T$. We must identify the timesteps which correspond to times we are confident that this simulation is still resolved. We assume the transfer of energy is largely local. This means that the energy that begins in low-frequency modes at the beginning will begin to drain into modes with increasing $|\mathbf{k}|$ as time evolves. The assumption of local energy transfer means that we can assume a simulation is still resolved as long as energy has not yet reached the edge of the computational domain.

We process the data to find a conservative estimate of the resolved timesteps. At each timestep, we define $\hat{\mathbf{u}}(t)$ to be a restricted version of the calculated solution of size $F = \{\mathbf{k}\;\;|\;\;\mathbf{k}\in[-M/2,M/2-1]^3\}$ at the observed timesteps. We can calculate to first order the rate of change of energy flowing out of these modes through the $t$-model:
\begin{equation}
\Delta E_F(t) = t \bigg[ \sum_{\mathbf{k}\in F}\mathbf{R}^1_{\mathbf{k}}(\hat{\mathbf{u}})\cdot\overline{\mathbf{u}}_{\mathbf{k}}+\mathbf{u}_{\mathbf{k}}\cdot\overline{\mathbf{R}}^1_{\mathbf{k}}(\hat{\mathbf{u}}) \bigg].
\end{equation}We restrict ourselves to timesteps where the full data are still resolved by limiting ourselves to a specific set of timesteps $t^* = \{t\;\;|\;\;10^{-16}<\Delta E_F(t) < 10^{-10}\}$. If there is so little energy leaving the cube corresponding to the inner 1/27th of the simulated domain, we conclude that the energy has not yet cascaded to the edge of the domain and the simulation is still resolved. We double-checked this by computing:
\begin{equation}
\Delta E_G(t) = t \bigg[ \sum_{\mathbf{k}\in G}\mathbf{R}^1_{\mathbf{k}}(\mathbf{u})\cdot\overline{\mathbf{u}}_{\mathbf{k}}+\mathbf{u}_{\mathbf{k}}\cdot\overline{\mathbf{R}}^1_{\mathbf{k}}(\mathbf{u}) \bigg],
\end{equation}which is a first-order approximation of the amount of energy that \emph{should} be flowing out of the full system were it larger. We found that for all $t\in t^*$ this quantity was below machine precision. This makes us more confident that the timesteps $t^*$ can be trusted as resolved. Thus, we can use it for fitting our reduced order models.

Consider a reduced order model of resolution $N$ that includes reduced order models up through order $n$. By this we mean the ROM would compute the solutions for the wavenumbers in $F_N = \{\mathbf{k} \;\;|\;\;\mathbf{k}\in [-N,N-1]^3\}$. We assume the renormalization coefficients depend upon the system size $N$. We computed these coefficients using a least squares fit. For the algebraically decaying coefficients $\alpha_i(t)=a_it^{-i}$, we minimized:
\begin{equation}
C_{N,n}(\mathbf{a}) = \sum_{\mathbf{k}\in F_N}\sum_{t\in t^*} \left(\Delta E_{\mathbf{k}} - \Delta E_{\mathbf{k}}^0-\sum_{i=1}^{n}a_i \Delta E_{\mathbf{k}}^i\right)^2\label{KdV_scaling_cost}
\end{equation}where $\mathbf{a}$ is the vector of renormalization coefficients. For the constant renormalization coefficients $\alpha_i'(t) = a_i'$, we minimized:
\begin{equation}
C'_{N,n}(\mathbf{a}') = \sum_{\mathbf{k}\in F_N}\sum_{t\in t^*} \left(\Delta E_{\mathbf{k}} - \Delta E_{\mathbf{k}}^0- \sum_{i=1}^{n}a'_it^i \Delta E_{\mathbf{k}}^i\right)^2.\label{Burgers_scaling_cost}
\end{equation}

\begin{figure}[h]
\begin{center}
\includegraphics[width=0.45\textwidth]{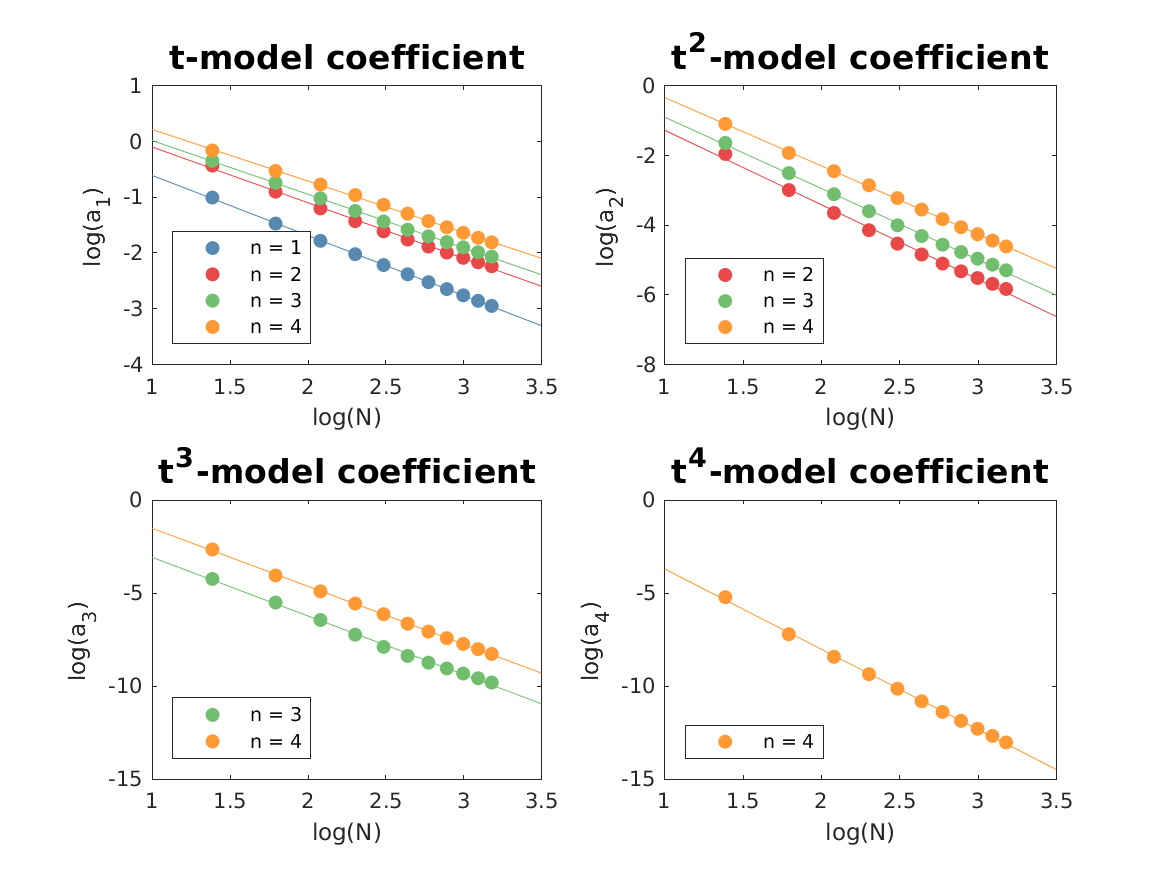}  \includegraphics[width=0.45\textwidth]{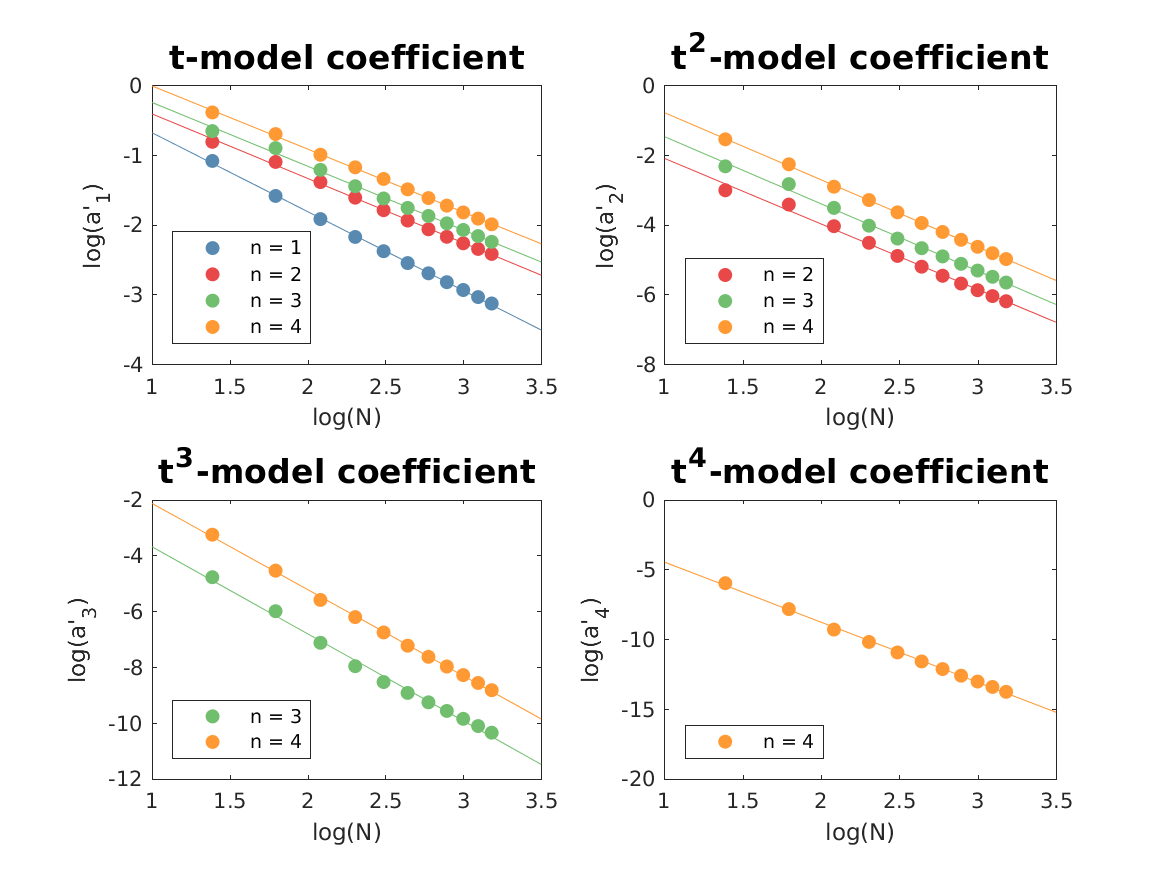}
\end{center}
\caption{Optimal algebraically decaying renormalization coefficients $a_i$ (left) and constant renormalization coefficients $a_i'$ (right) plotted on a log-log scale against the system resolution $N$. The optimal coefficients were calculated by minimizing \eqref{KdV_scaling_cost} and \eqref{Burgers_scaling_cost}, respectively, for $N = 4,6,\dots,24$ and $n = 1,2,3,4$ against data produced by a still-resolved $M=48$ full simulation.}\label{fig:Euler_scaling}
\end{figure}

We ran a full simulation of size $M = 48$ and used it to compute renormalization coefficients for reduced order models of size $N = 4,6,\dots,24$ with algebraically decaying renormalization coefficients and constant renormalization coefficients. In each case, we considered models that included up through $n = 1$,$2$,$3$,$4$ terms from the complete memory approximation. In all cases, there was an apparent algebraic dependence upon $N$ for the constants $\mathbf{a}$ and $\mathbf{a}'$ as depicted in Fig. \ref{fig:Euler_scaling}. We also used simulations of size $M = 24$ and $M = 32$ to compute renormalization coefficients for reduced order models of size $N = 4,6,\dots,12$ and $N = 4,6,\dots,16$, respectively. The coefficients found in each of these calculations were similar. As the largest simulation, we trust the results from the $M = 48$ simulation the most.

Due to the apparent algebraic dependence, we conclude that the functional form of the renormalization coefficients is:
\begin{align}
\alpha_i(t) &= \beta_i^nN^{\gamma_i^n} t^{-i}\\
\alpha_i'(t) &= {\beta}_i^{\prime n}{N}^{\gamma_i^{\prime n}}.
\end{align}We identified the parameters $\beta_i^n$, $\gamma_i^n$, ${\beta}_i^{\prime n}$ and $\gamma_i^{\prime n}$ by computing a linear least squares fit of $\log(\alpha_i)$ and $\log(N)$. The resulting coefficients and the correlation coefficient of the linear least-squares fit $r^2$ are presented in Table \ref{tab:scaling_law_table}. 

There are a few comments we can make from these data. First, we see from the correlation coefficients that these fits are quite good. It appears the scaling law form is a good representation of the functional form of the renormalization coefficients. We see also that the exponent in those scaling laws ($\gamma_i^n$ and $\gamma_i^{\prime n}$) seem relatively independent of the number of terms $n$ included in the reduced model. The prefactors $\beta_i^n$ and $\beta_i^{\prime n}$, however, slowly grow in magnitude as more terms are included. We also see that the coefficients of the even terms are negative while the coefficients of the odd terms are positive in all cases. Thus, the renormalized coefficients agree with the unrenormalized coefficients in sign though not in magnitude.

\begin{table}[H]
\begin{center}
\begin{adjustbox}{max width=\textwidth}
\begin{tabular}{|c||c|c|c|c||c|c|c|c||c|c|c|c|}
\hline$n$ & $\beta_1^n$& $\beta_2^n$& $\beta_3^n$& $\beta_4^n$& $\gamma_1^n$& $\gamma_2^n$& $\gamma_3^n$& $\gamma_4^n$ & $r_1$& $r_2$ & $r_3$ & $r_4$\\
\hhline{|=#=|=|=|=#=|=|=|=#=|=|=|=|}
1 & $1.591$ &  &  &  & $-1.077$ &  &  &  & $1.000$ &  &  &  \\
\hline
2 & $2.448$ & $-2.341$ &  &  & $-0.999$ & $-2.136$ &  &  & $0.998$ & $0.996$ &  &  \\
\hline
3 & $2.650$ & $-3.094$ & $1.068$ &  & $-0.962$ & $-2.042$ & $-3.148$ &  & $1.000$ & $0.999$ & 0.$998$ &  \\
\hline
4 & $3.110$ & $-5.006$ & $4.924$ & $-1.828$ & $-0.924$ & $-1.959$ & $-3.115$ & $-4.312$ & $0.999$ & $1.000$ & $1.000$ & $0.999$ \\
\hline
\end{tabular}
\end{adjustbox}

\vspace{1.5em}

\begin{adjustbox}{max width=\textwidth}
\begin{tabular}{|c||c|c|c|c||c|c|c|c||c|c|c|c|}
\hline$n$ & $\beta_1^{\prime n}$& $\beta_2^{\prime n}$& $\beta_3^{\prime n}$& $\beta_4^{\prime n}$& $\gamma_1^{\prime n}$& $\gamma_2^{\prime n}$& $\gamma_3^{\prime n}$& $\gamma_4^{\prime n}$ & $r_1'$& $r_2'$ & $r_3'$ & $r_4'$\\
\hhline{|=#=|=|=|=#=|=|=|=#=|=|=|=|}
1 & $1.574$ &  &  &  & $-1.132$ &  &  &  & $0.999$ &  &  &  \\
\hline
2 & $1.677$ & $-0.805$ &  &  & $-0.926$ & $-1.879$ &  &  & $0.998$ & $0.994$ &  &  \\
\hline
3 & $1.955$ & $-1.570$ & $0.573$ &  & $-0.915$ & $-1.925$ & $-3.122$ &  & $0.997$ & $0.996$ & $0.994$ &  \\
\hline
4 & $2.454$ & $-3.124$ & $2.628$ & $-0.860$ & $-0.906$ & $-1.924$ & $-3.091$ & $-4.306$ & $0.999$ & $1.000$ & $0.999$ & $0.998$ \\
\hline
\end{tabular}
\end{adjustbox}
\caption{Scaling laws that approximate the observed optimal correlation coefficients for 3D Euler's equations. The top table contains scaling laws for algebraically decaying renormalization coefficients $\alpha_i(t) = \beta_i^nN^{\gamma_i^n} t^{-i}$ while the bottom table contains scaling laws for constant renormalization coefficients $\alpha_i'(t) = {\beta}_i^{\prime n}{N}^{\gamma_i^{\prime n}}$. We computed an $M = 48$ full simulation to construct $\Delta E_{\mathbf{k}}$. The algebraically decaying renormalization coefficient scaling laws were found by minimizing \eqref{KdV_scaling_cost} with $n = 1,2,3,4$ and $N = 4,6,\dots,24$, then conducting a linear least squares fit of $\log(a_i)$ against $\log(N)$. The correlation coefficient $r^2$ of this log-log fit is also provided. The same data methodology was used to compute the constant renormalization coefficient scaling laws.} \label{tab:scaling_law_table}
\end{center}
\end{table}

\section{Results}

The behavior of the solution to the three-dimensional Euler's equations with a smooth initial condition remains unknown. Consequently, we cannot compare the results of our ROMs to the exact solution for accuracy. Instead, we endeavour to produce ROMs that remain \emph{stable} over a long time. We will have to rely upon secondary means of inferring the accuracy of the resultant ROMs. Our results, not fully validated as they are, can be interpreted as evidence that is suggestive of long-term behavior of a subset of Fourier modes evolved according to Euler's equations.

\begin{figure}[h]
\begin{center}
 \includegraphics[width=0.5\textwidth]{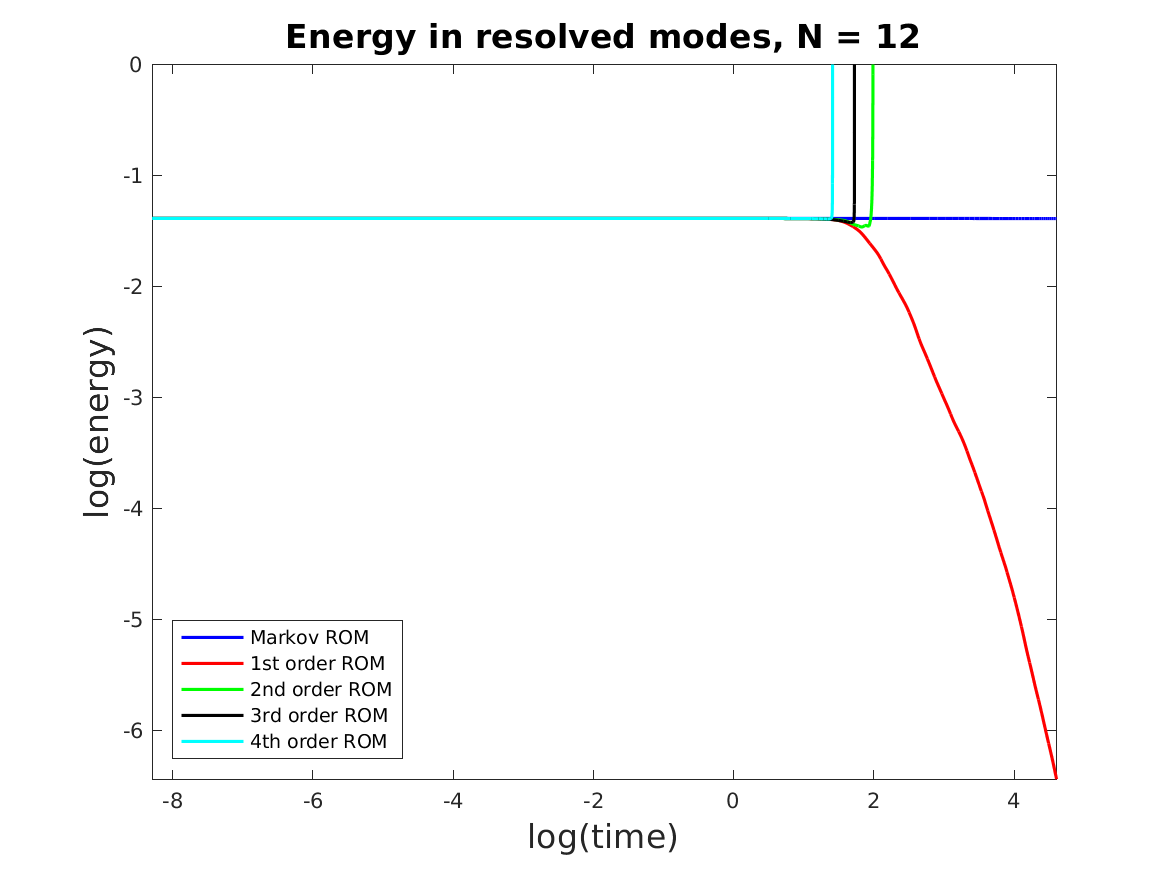}
\end{center}
\caption{The energy contained in the resolved modes of several ROMs of resolution $N=12$ using \emph{constant} renormalization coefficients as described in the bottom of Table \ref{tab:scaling_law_table} depicted on a log-log plot. The Markov model is stable but does not drain any energy. The first order ROM is stable by construction. All other ROMs are unstable, and the time of instability grows smaller as the order of the ROM increases. These results are qualitatively the same for other resolutions $N$.}\label{fig:unstable_burgers_ROM}
\end{figure}

First, it should be noted that renormalized ROMs with \emph{constant} renormalization coefficients proved to be unstable for $n>1$. The $t$-model is stable, by construction, but the addition of higher-order terms, even when renormalized, rendered the simulations unstable. As the order of the ROM increases for a fixed resolution, the time at which the model becomes unstable becomes earlier (see Fig. \ref{fig:unstable_burgers_ROM}). Consequently, we conclude that the constant renormalization coefficients are not the correct choice for producing stable ROMs for Euler's equations.

\begin{figure}[h]
\begin{center}
 \includegraphics[width=0.75\textwidth]{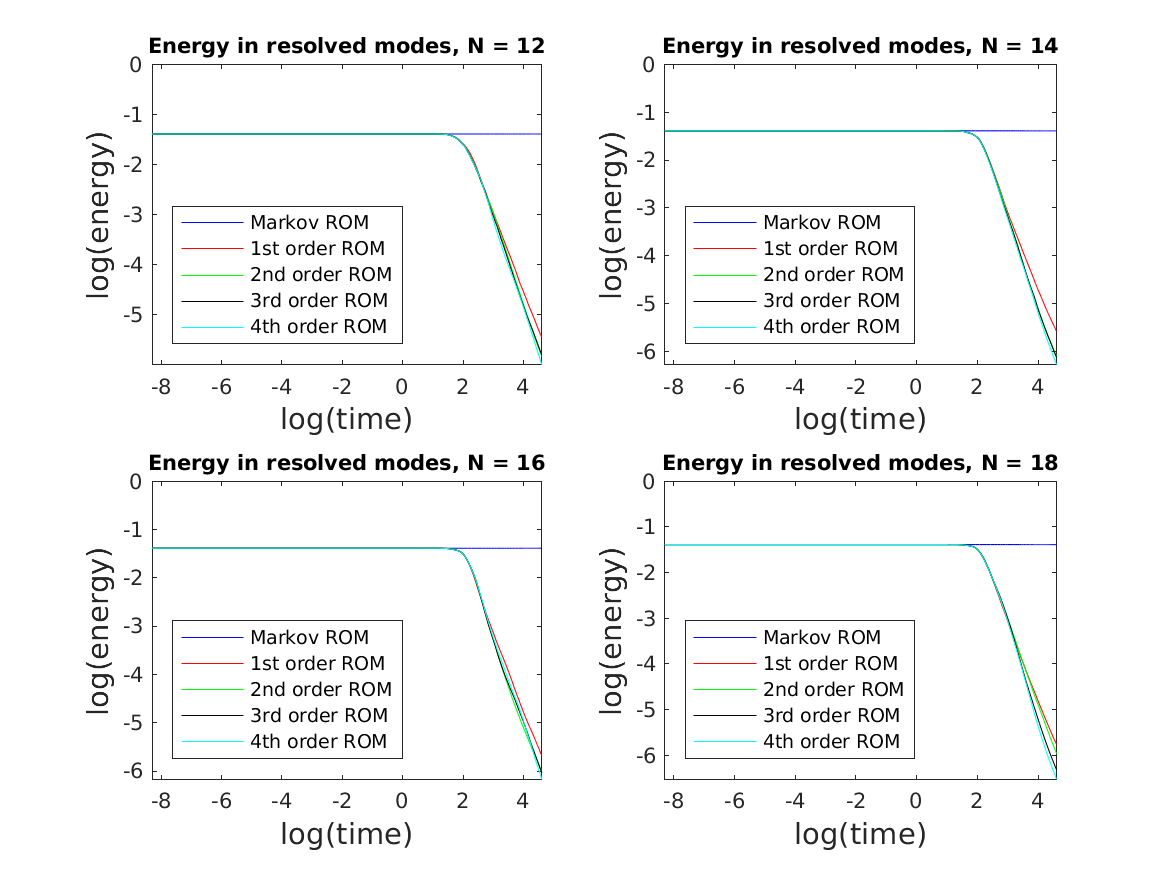}
\end{center}
\caption{The energy contained in the resolved modes several ROMs of resolution $N=12,14,16,18$ up to time $t=100$ depicted on a log-log plot. The Markov model does not drain energy at all. The other four ROMs use algebraically decaying renormalization coefficients as described in Table \ref{tab:scaling_law_table}. Note that the behavior of the energy appears to converge as the order of the model increases, indicating that we are in the perturbative regime. The results for other resolutions $N$ are qualitatively similar and the convergence appears to be faster as $N$ grows larger.}\label{fig:perturbative_euler}
\end{figure}

On the other hand, the ROMs with renormalization coefficients that \emph{decay algebraically with time} led to solutions that remained stable until at least $t=1000$. When we fix a resolution $N$ and simulate ROMs that include up through degree $n=1,2,3,4$, as seen in Fig. \ref{fig:perturbative_euler}, the results converge quickly with increasing order, which suggests that we are indeed in a perturbative regime. Each additional term in a ROM is more expensive to compute, and the fast convergence gives us confidence that including additional terms will only minimally affect our results. Thus, we will assume that the fourth order ROMs represent the most accurate simulations of the dynamics of the resolved modes.

\begin{figure}[h]
\begin{center}
 \includegraphics[width=0.5\textwidth]{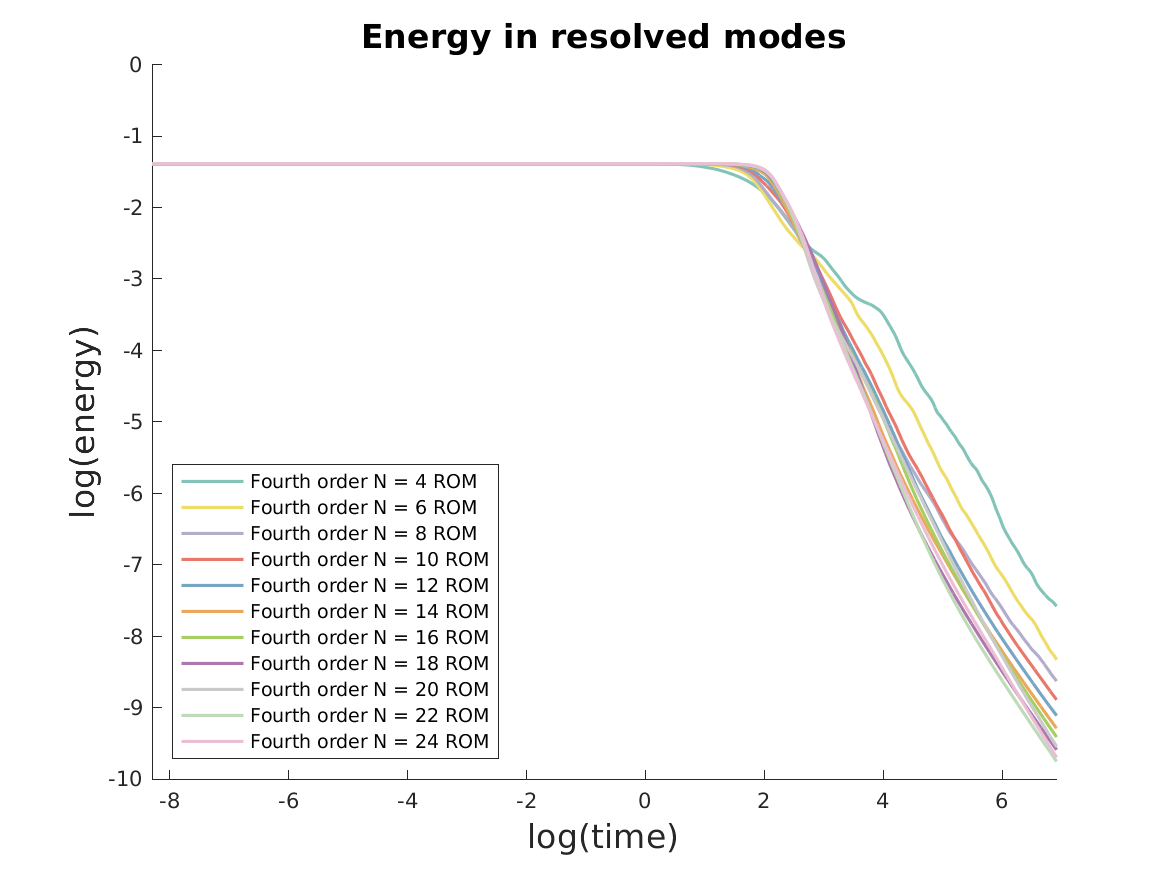}
\end{center}
\caption{The energy contained in the resolved modes of fourth order ROMs of resolutions $N=4,6,\dots,24$ up to time $t=1000$ depicted on a log-log plot. These simulations use algebraically decaying renormalization coefficients as described in Table \ref{tab:scaling_law_table}. The qualitative behavior is the same for each resolution: energy does not drain at first, until it suddenly begins to do so at a nearly constant rate. After a long time, it gradually shifts to a different constant drain rate. The time at which the energy drain begins becomes larger as $N$ increases, and the initial slope grows steeper with increasing $N$. Both appear to have limiting values. The total energy ejected by the end of the simulation grows monotonically with increasing $N$.}\label{fig:long_energy}
\end{figure}

As one increases the resolution of our fourth order ROMs, several fascinating patterns emerge. First, we consider the energy contained in the resolved modes. The Markov term conserves energy, so any draining of energy is accomplished by the memory terms alone. Fig. \ref{fig:long_energy} depicts the energy decay of ROMs with resolution $N=4,6,\dots,24$ up to time $t=1000$ on a log-log plot. We see that in all cases there is monotonic energy decay. As time goes on, the results become stratified: the amount of energy remaining in the system \emph{decreases} as the resolution of the model grows. This indicates significant activity in the high-frequency modes that increases with the resolution. This is one point of evidence for a singularity. Were there not a singularity, one would not expect each larger ROM to drain more energy. 

We also see several other patterns as the resolution increases.  The decay of energy is relatively linear on a log-log plot, indicating algebraic energy ejection from the resolved modes. The slope indicates the exponent of the decay. We computed the slope from the data for which between 50\% and 90\% of the initial energy has left the system. These slopes seem to converge towards approximately -2.5. Finally, we see that the rate of energy ejection eventually becomes slightly smaller. We computed the slope from the data after 99.5\% of the initial energy had left the system. These new ejection rates seem to converge towards -1.5. All these observations are enumerated in Table \ref{tab:energy_decay_Euler}.

\begin{table}[H]
\begin{center}
\begin{tabular}{|c|c|c|}
\hline ROM resolution $N$ &  Initial decay rate& Second decay rate\\
\hline
4 &  -0.808 & -1.208\\
\hline
6 &  -1.003 & -1.291\\
\hline
8 &  -1.575 & -1.163\\
\hline
10 &  -1.646 & -1.298\\
\hline
12 &  -1.916 & -1.274\\
\hline
14 &  -1.981 & -1.265\\
\hline
16 &  -2.005 & -1.357\\
\hline
18 &  -1.989 & -1.298\\
\hline
20 &  -2.208 & -1.485\\
\hline
22 & -2.226 & -1.348\\
\hline
24 &  -2.414 & -1.413\\
\hline
\end{tabular}
\caption{Observations of the energy decay in reduced order models of 3D Euler's equations. Simulations use algebraically decaying renormalization coefficients as described in the top of Table \ref{tab:scaling_law_table}. When plotted on a log-log plot, the energy begins to decay at a fixed rate. The initial decay rate is the slope of a least-squares fit line to the log-log data for which between 50\% and 90\% of the initial energy has left the resolved modes. The second decay rate is the slope of a least-squares fit line to the log-log data after 99.5\% of the data has left the resolved modes.}\label{tab:energy_decay_Euler}
\end{center}
\end{table}

Next, we computed several other interesting dynamic quantities from our ROMs. In each of these cases, we found perturbative convergence as the order of the model increases. This suggests that we can continue to treat the fourth order ROM as the most accurate result. 
\begin{figure}[h]
\begin{center}
 \includegraphics[width=0.5\textwidth]{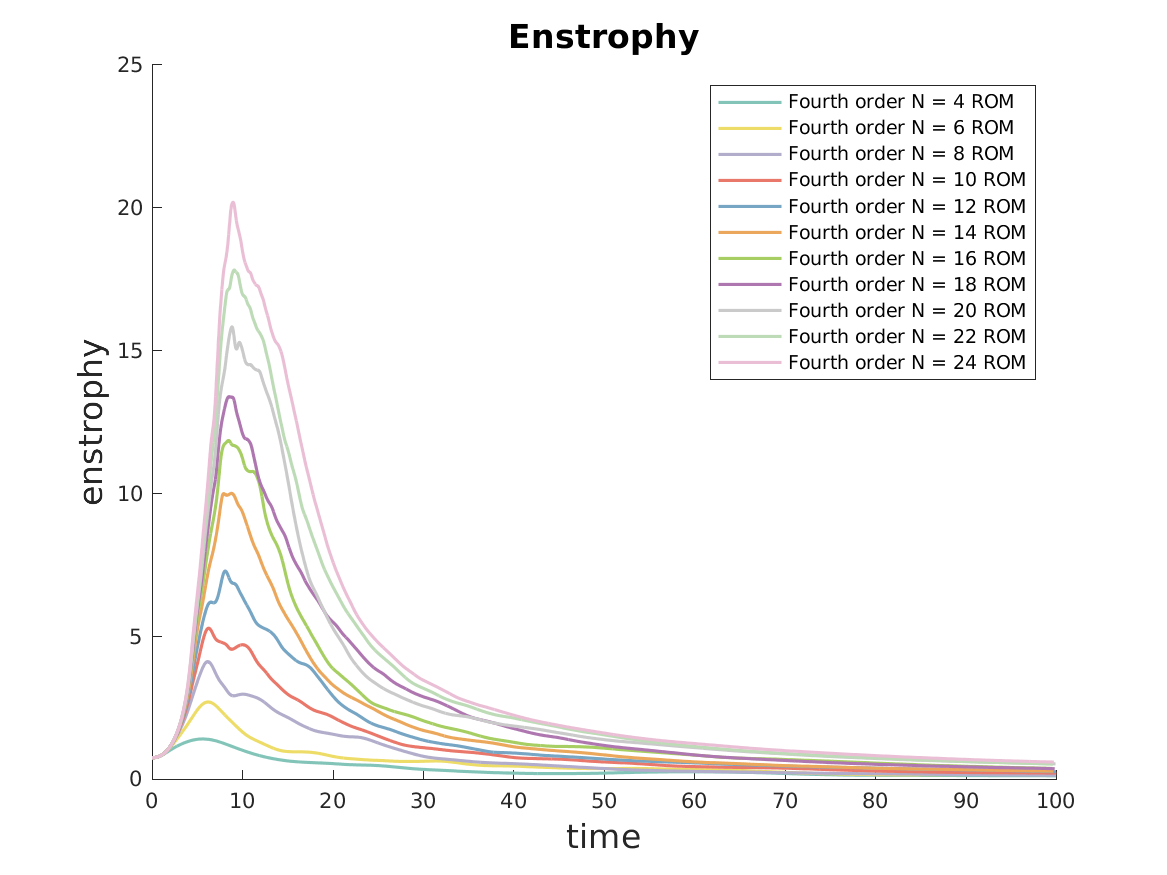}
\end{center}
\caption{The enstrophy plotted against time for fourth order ROMs of size $N=4,6,\dots,24$ using algebraically decaying renormalization coefficients as described in Table \ref{tab:scaling_law_table} up to time $t=100$. The enstrophy for each ROM begins small, grows to a maximum at a finite time, and then decays. As $N$ increases, the maximum value achieved increases, and the time at which this maximum is achieved appears to coincide with the time at which energy begins draining from the system. If this pattern continues and the enstrophy approaches infinity at some finite time as the number of simulated modes increases, we can conclude that a finite-time singularity occurs.}\label{fig:enstrophy}
\end{figure}
\begin{figure}[h]
\begin{center}
 \includegraphics[width=0.5\textwidth]{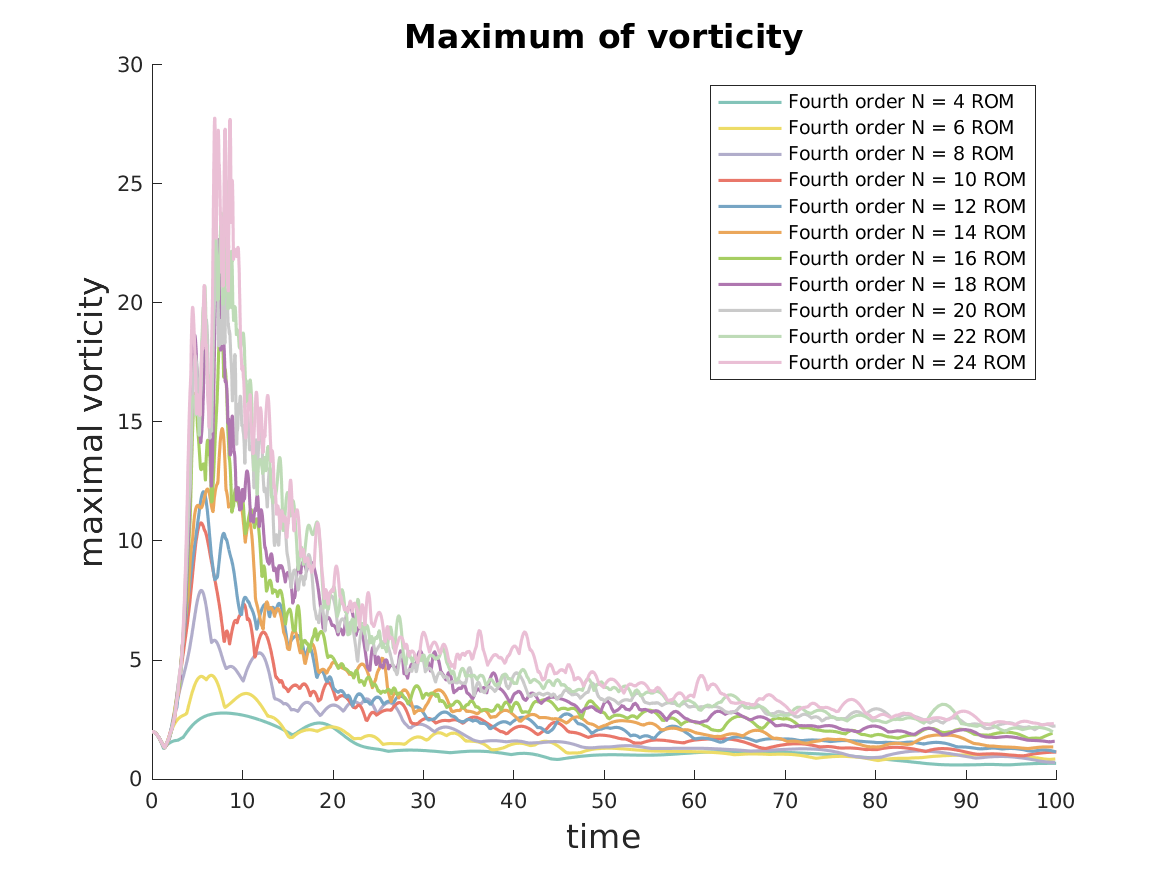}
\end{center}
\caption{The maximal vorticity plotted against time for fourth order ROMs of size $N=4,6,\dots,24$ using algebraically decaying renormalization coefficients as described in Table \ref{tab:scaling_law_table} up to time $t=100$. The maximal vorticity grows to a maximum at a finite time like the enstrophy. As the resolution increases, the maximum value achieved increases. This is further suggestive that a finite-time singularity occurs.}\label{fig:vorticity}
\end{figure}

\begin{figure}[h]
\begin{center}
 \includegraphics[width=0.5\textwidth]{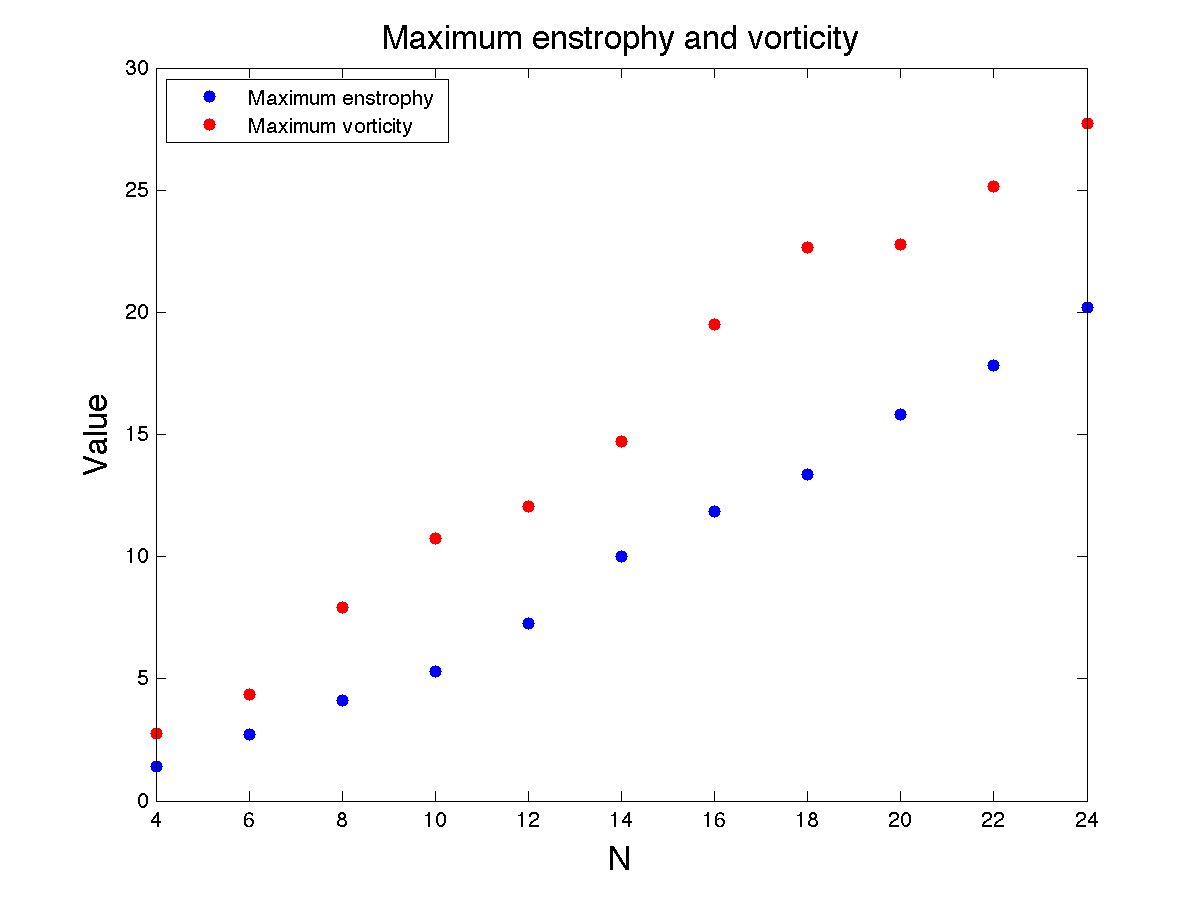}
\end{center}
\caption{The maximum (over time) of the enstrophy and of the maximal vorticity as a function of the resolution. As the resolution increases, the maxima increase almost linearly.}\label{fig:maxima}
\end{figure}

The curl of the solution is called the \emph{vorticity} in fluid dynamics. The enstrophy is defined as:
\begin{equation}
e(t) = \int |\nabla\times\mathbf{u}|^2\,\mathrm{d}\mathbf{x},
\end{equation}
where $\nabla\times\mathbf{u}$ is the vorticity. If a singularity does occur, we would expect a peak in the enstrophy to occur at the time of the singularity. Furthermore, we would expect this peak to grow larger as the resolution is increased (and become infinite as the resolution grows to infinity). Indeed, these are exactly the results we observe (Fig. \ref{fig:enstrophy}). Furthermore, we observe that the time at which the enstrophy peaks roughly coincides with the time at which the energy begins to flow out of the system at a fixed algebraic rate. The maximum of the vorticity
\begin{equation} ||\mathbf{\omega}||_{\infty}(t) = \max|\nabla\times \mathbf{u}|
\end{equation} is the best indicator of singular behavior (see e.g. \cite{doering1995applied}). During a singularity this quantity will blow up even if the enstrophy does not. In the presence of a singularity, one would again expect a peak in the maximum of the vorticity at some finite time, and we would expect this peak to grow towards infinity as the resolution is increased. We observe this here as well, and again the time of the peak seems to roughly coincide with the point at which energy begins flowing out of the system (Fig. \ref{fig:vorticity}). In Fig. \ref{fig:maxima} we plot the maximal values (over time) of the enstrophy and of the maximum vorticity as a function of resolution. These maxima increase almost linearly with resolution. This trend is consistent with  the formation of a finite-time singularity.

\section{Discussion}

The exact behavior of solutions to 3D Euler's equation is unknown. Even modern simulations with exceptionally high resolution cannot proceed for long times. Thus, our reduced order models represent an advancement in the ability to simulate these equations. Without an exact solution to validate against, it is difficult to ascertain whether our results are accurate in addition to stable. However, there are a few hints: the convergence of behavior in Fig. \ref{fig:perturbative_euler} indicates that our ROMs have a perturbative structure. That is, each additional order in the ROM modifies the solution less and less. It appears to be converging towards \emph{something} and it is not unreasonable to think that it is the exact solution. Next, Table \ref{tab:scaling_law_table} demonstrates that adding additional terms does not significantly change the scaling laws for the previous terms. Each additional term is making \emph{corrections} to previously captured behavior, but their contributions seem to be somewhat orthogonal to one another. This behavior gives us confidence in trusting our results.

The simulations are rather expensive. Each convolution requires a three-dimensional FFT and IFFT. Furthermore, for a system of resolution $N$, the FFTs are of size $(2\times 2\times 3/2 \times N)^3$.  One factor of 2 comes from the fact that we include both positive and negative modes. The other factor of 2 because the terms in the reduced model require a full system double the size. The factor of $3/2$ is needed to dealias the results. The Markov model requires only one FFT. The first order ROM requires an additional two convolutions. The second order ROM adds another five, while the third order ROM has an additional nine on top of that. The fourth order ROM requires all past convolutions, plus another 20. Thus, the simulation cost grows very quickly as the degree of the ROM increases. On the other hand, as stated above, we are \emph{incapable} of performing a brute force calculation beyond a few units of time even on modern cutting-edge high-resolution simulations. With our models we can integrate out to time $t=1000$ for $N = 24$ in only a few days on a laptop computer. The results in Figs. \ref{fig:long_energy}, \ref{fig:enstrophy}, \ref{fig:vorticity}, and \ref{fig:maxima} and Table \ref{tab:energy_decay_Euler} all provide strong evidence of a singularity occurring in Euler's equations with the Taylor-Green initial condition. Our simulations are able to pass through the singularity, which would not be possible for a non-reduced system even with infinite computational power (if a singularity does indeed occur).

We have found that algebraically decaying renormalization coefficients are necessary to produce stable reduced order models. We found an algebraic dependence of the prefactors upon the resolution $N$. It is possible that a different time dependence for the renormalization coefficients could yield stable and accurate simulations. This is an area of future inquiry by the authors. The algebraic decay of energy from the resolved modes is suggestive of a finite time singularity. The increasing peaks of the enstrophy and maximum of the vorticity provide further evidence for the formation of a singularity. The convergence of each of these events to a fixed time leads us to conclude that these reduced order models provide strong evidence for a finite time singularity developing from a smooth initial condition in Euler's equations. This is an important and long-standing open problem in fluid dynamics. Reduced order models allow us to utilize the multiscale structure of problems to evolve only a subset of the variables in the system.

Our evidence is \emph{suggestive} of a finite-time singularity, but this evidence would become more persuasive with larger-scale simulations. We were only able to simulate relatively small $24^3$ mode simulations using our memory approximation. This precludes interesting visualization of the results because they lack sufficient resolution. Furthermore, it is possible that new and exciting effects only become apparent at sufficient resolution. We are interested in applying our reduced order models to a more finely-resolved simulation on a powerful computer. Note that to simulate the reduced model for a higher resolution we can use the estimated scaling laws to compute the renormalized coefficients. The ability to calibrate a reduced model using smaller and well-thought calculations and then extrapolate to hitherto unreachable regimes is a major goal of scientific computing. Finally, it will be very interesting to investigate how adding viscous dissipation (Navier-Stokes equations) will alter the scaling dependence of renormalized coefficients, including the likely occurrence of incomplete similarity \cite{barenblatt,chorin2005spectra,price2017renormalized}.

\section{Acknowledgements}

The work of PS was supported by the U.S. Department of Energy (DOE) Office of Science, Office of Advanced Scientific Computing Research (ASCR) as part of the Multifaceted Mathematics for Rare, Extreme Events in Complex Energy and Environment Systems (MACSER) project. Pacific Northwest National Laboratory is operated by Battelle for the DOE under Contract DE-AC05-76RL01830.

\bibliographystyle{plain}
\bibliography{MultiScale}

\begin{thebibliography}{10}

\bibitem{agafontsev2015development}
DS~Agafontsev, EA~Kuznetsov, and AA~Mailybaev.
\newblock Development of high vorticity structures in incompressible 3d {E}uler
  equations.
\newblock {\em Physics of Fluids}, 27(8):085102, 2015.

\bibitem{ayala2017extreme}
Diego Ayala and Bartosz Protas.
\newblock Extreme vortex states and the growth of enstrophy in
  three-dimensional incompressible flows.
\newblock {\em Journal of Fluid Mechanics}, 818:772--806, 2017.

\bibitem{barenblatt}
Grigory~I Barenblatt.
\newblock {\em Scaling}.
\newblock Cambridge University Press, 2003.

\bibitem{bernstein2007optimal}
David Bernstein.
\newblock Optimal prediction of {B}urgers's equation.
\newblock {\em Multiscale Modeling \& Simulation}, 6(1):27--52, 2007.

\bibitem{bustamante2012euler}
Miguel Bustamante and Marc Brachet.
\newblock Interplay between the {B}eale-{K}ato-{M}ajda theorem and the
  analyticity-strip method to investigate numerically the incompressible
  {E}uler singularity problem.
\newblock {\em Physical Review E}, 86:066302, 2012.

\bibitem{chandy2009t}
Abhilash~J Chandy and Steven~H Frankel.
\newblock The t-model as a large eddy simulation model for the
  {N}avier-{S}tokes equations.
\newblock {\em Multiscale Modeling \& Simulation}, 8(2):445--462, 2009.

\bibitem{C94}
Alexandre~J Chorin.
\newblock {\em Vorticity and turbulence}, volume 103.
\newblock Springer Science \& Business Media, 1994.

\bibitem{chorin2005spectra}
Alexandre~J Chorin and Ole~H Hald.
\newblock Viscosity-dependent inertial spectra of the {B}urgers and
  {K}orteweg-de{V}ries-{B}urgers equations.
\newblock {\em Proceedings of the National Academy of Sciences},
  102(11):3921--3923, 2005.

\bibitem{chorin2000optimal}
Alexandre~J Chorin, Ole~H Hald, and Raz Kupferman.
\newblock Optimal prediction and the {M}ori-{Z}wanzig representation of
  irreversible processes.
\newblock {\em Proceedings of the National Academy of Sciences},
  97(7):2968--2973, 2000.

\bibitem{chorin2002optimal}
Alexandre~J Chorin, Ole~H Hald, and Raz Kupferman.
\newblock Optimal prediction with memory.
\newblock {\em Physica D: Nonlinear Phenomena}, 166(3):239--257, 2002.

\bibitem{chorin2007problem}
Alexandre~J Chorin and Panos Stinis.
\newblock Problem reduction, renormalization, and memory.
\newblock {\em Communications in Applied Mathematics and Computational
  Science}, 1(1):1--27, 2007.

\bibitem{constantin1994geometric}
Peter Constantin.
\newblock Geometric statistics in turbulence.
\newblock {\em SIAM Review}, 36(1):73--98, 1994.

\bibitem{constantin2008euler}
Peter Constantin.
\newblock {E}uler and {N}avier-{S}tokes equations.
\newblock {\em Publicacions Matem{\`a}tiques}, pages 235--265, 2008.

\bibitem{constantin2017analysis}
Peter Constantin.
\newblock {\em Analysis of Hydrodynamic Models}.
\newblock {SIAM}, 2017.

\bibitem{deng2006improved}
Jian Deng, Thomas~Y Hou, and Xinwei Yu.
\newblock Improved geometric conditions for non-blowup of the 3{D}
  incompressible {E}uler equation.
\newblock {\em Communications in Partial Differential Equations},
  31(2):293--306, 2006.

\bibitem{doering1995applied}
Charles~R Doering and John~D Gibbon.
\newblock {\em Applied analysis of the {N}avier-{S}tokes equations}, volume~12.
\newblock Cambridge University Press, 1995.

\bibitem{elgindi2018finite}
Tarek~M Elgindi and In-Jee Jeong.
\newblock Finite-time singularity formation for strong solutions to the
  axi-symmetric 3{D} {E}uler equations.
\newblock {\em arXiv preprint arXiv:1802.09936}, 2018.

\bibitem{georgi1993}
Howard Georgi.
\newblock Effective {F}ield {T}heory.
\newblock {\em Annual Review Nuclear Particle Science}, 43:209--252, 1993.

\bibitem{goldenfeld1992}
Nigel Goldenfeld.
\newblock {\em Lectures on {P}hase {T}ransitions and the {R}enormalization
  {G}roup}.
\newblock Perseus Books, 1992.

\bibitem{grafke2013lagrangian}
Tobias Grafke and Rainer Grauer.
\newblock Lagrangian and geometric analysis of finite-time {E}uler
  singularities.
\newblock {\em Procedia IUTAM}, 9:32--56, 2013.

\bibitem{hald2007optimal}
Ole~H Hald and Panos Stinis.
\newblock Optimal prediction and the rate of decay for solutions of the {E}uler
  equations in two and three dimensions.
\newblock {\em Proceedings of the National Academy of Sciences},
  104(16):6527--6532, 2007.

\bibitem{hou2006dynamic}
Thomas~Y Hou and Ruo Li.
\newblock Dynamic depletion of vortex stretching and non-blowup of the 3-{D}
  incompressible {E}uler equations.
\newblock {\em Journal of Nonlinear Science}, 16(6):639--664, 2006.

\bibitem{isett2017onsager}
Philip Isett.
\newblock {\em H\"older Continuous {E}uler Flows in Three Dimensions with
  Compact Support in Time}.
\newblock Princeton University Press, 2017.

\bibitem{koopman1931hamiltonian}
Bernard~O Koopman.
\newblock {H}amiltonian systems and transformation in {H}ilbert space.
\newblock {\em Proceedings of the National Academy of Sciences},
  17(5):315--318, 1931.

\bibitem{lei2016}
H~Lei, NA~Baker, and X~Li.
\newblock Data-driven parameterization of the generalized {L}angevin equation.
\newblock {\em Proceedings of the National Academy of Sciences},
  33(3):14183--14188, 2016.

\bibitem{li2008blow}
Dong Li and Yakov~G Sinai.
\newblock Blow ups of complex solutions of the 3d {N}avier-{S}tokes system and
  renormalization group method.
\newblock {\em Journal of the European Mathematical Society}, 10(2):267--313,
  2008.

\bibitem{li2015memory}
Zhen Li, Xin Bian, Xiantao Li, and George~Em Karniadakis.
\newblock Incorporation of memory effects in coarse-grained modeling via the
  {M}ori-{Z}wanzig formalism.
\newblock {\em The Journal of Chemical Physics}, 143(24):3128, 2015.

\bibitem{luo2014potentially}
Guo Luo and Thomas~Y Hou.
\newblock Potentially singular solutions of the 3d axisymmetric {E}uler
  equations.
\newblock {\em Proceedings of the National Academy of Sciences},
  111(36):12968--12973, 2014.

\bibitem{majda}
Andrew~J Majda and Andrea~L Bertozzi.
\newblock {\em Vorticity and incompressible flow}, volume~27.
\newblock Cambridge University Press, 2002.

\bibitem{marchioro1994euler}
Carlo Marchioro and Mario Pulvirenti.
\newblock {\em Mathematical theory of incompressible nonviscous fluids},
  volume~96.
\newblock Springer, 1994.

\bibitem{moffatt}
HK~Moffatt, S~Kida, and Koji Ohkitani.
\newblock Stretched vortices--the sinews of turbulence; large-{R}eynolds-number
  asymptotics.
\newblock {\em Journal of Fluid Mechanics}, 259:241--264, 1994.

\bibitem{orlandi2012vortex}
P~Orlandi, S~Pirozzoli, and GF~Carnevale.
\newblock Vortex events in {E}uler and {N}avier-{S}tokes simulations with
  smooth initial conditions.
\newblock {\em Journal of Fluid Mechanics}, 690:288--320, 2012.

\bibitem{parish2017}
EJ~Parish and K~Duraisamy.
\newblock Non-{M}arkovian closure models for large eddy simulations using the
  {M}ori-{Z}wanzig formalism.
\newblock {\em Physical Review Fluids}, 2(1):014604, 2017.

\bibitem{MZrepos}
Jacob Price.
\newblock Dissertation release of renormalized {M}ori-{Z}wanzig git repository.
\newblock \url{http://doi.org/10.5281/zenodo.1246871}, 2018.

\bibitem{price2017renormalized}
Jacob Price and Panos Stinis.
\newblock Renormalized reduced order models with memory for long time
  prediction.
\newblock {\em arXiv preprint arXiv:1707.01955}, 2017.

\bibitem{shu2005euler}
Chi-Wang Shu, Wai-Sun Don, David Gottlieb, Oleg Schilling, and Leland Jameson.
\newblock Numerical convergence study of nearly incompressible, inviscid
  {T}aylor-{G}reen vortex flow.
\newblock {\em SIAM Journal on Scientific Computing}, 24(1):1--27, 2005.

\bibitem{stinis2007higher}
Panos Stinis.
\newblock Higher order {M}ori-{Z}wanzig models for the {E}uler equations.
\newblock {\em Multiscale Modeling \& Simulation}, 6(3):741--760, 2007.

\bibitem{stinis2012numerical}
Panos Stinis.
\newblock Numerical computation of solutions of the critical nonlinear
  {S}chr{\"o}dinger equation after the singularity.
\newblock {\em Multiscale Modeling \& Simulation}, 10(1):48--60, 2012.

\bibitem{stinis2013renormalized}
Panos Stinis.
\newblock Renormalized reduced models for singular {PDE}s.
\newblock {\em Communications in Applied Mathematics and Computational
  Science}, 8(1):39--66, 2013.

\bibitem{stinis2015renormalized}
Panos Stinis.
\newblock Renormalized {M}ori-{Z}wanzig-reduced models for systems without
  scale separation.
\newblock In {\em Proceedings of the Royal Society of London A: {M}athematical,
  Physical and Engineering Sciences}, volume 471, page 20140446, 2015.

\bibitem{tao2016blow}
Terence Tao.
\newblock Finite time blowup for {L}agrangian modifications of the
  three-dimensional {E}uler equation.
\newblock {\em Annals of PDE}, 2(9), 2016.

\bibitem{zwanzig1961memory}
Robert Zwanzig.
\newblock Memory effects in irreversible thermodynamics.
\newblock {\em Physical Review}, 124(4):983, 1961.

\end{thebibliography}

\newpage

\appendix 

\section{Derivation of higher-ordered terms}\label{sec:higherorder}

We derived the functional form of the first and second order terms in the complete memory approximation of Euler's equations in the main text. Here, we demonstrate the derivation of the third order term and present the result for the fourth order term as produced by a symbolic computation.

\subsection{Third-order Term}

We will make use of the terms we have already computed in order to simplify our derivation of the third order term. Under the complete memory approximation, the third term is:
\[\mathbf{R}^3_{\mathbf{k}}(\hat{\mathbf{u}})  =  Pe^{t\mathcal{L}}P\mathcal{L}[PLPL-2PLQL-2QLPL+QLQL]Q\mathcal{L}\mathbf{u}_{\mathbf{k}}^0.
\]We rewrite it with no $Q\mathcal{L}$ terms.
\begin{align*}
\mathbf{R}^3_{\mathbf{k}}(\hat{\mathbf{u}})  =&  Pe^{t\mathcal{L}}P\mathcal{L}[PLPL-2PLQL-2QLPL+QLQL]Q\mathcal{L}\mathbf{u}_{\mathbf{k}}^0\\
=&Pe^{t\mathcal{L}}P\mathcal{L}[3P\mathcal{L}(2P\mathcal{L}-P\mathcal{L})-3\mathcal{L}P\mathcal{L}+\mathcal{L}\mathcal{L}]Q\mathcal{L}\mathbf{u}_{\mathbf{k}}^0
\end{align*}We recognize that we have already computed an expression for $P\mathcal{L}(2P\mathcal{L}-\mathcal{L})Q\mathcal{L}\mathbf{u}_{\mathbf{k}}^0$ during our derivation of the second order term. This leaves two additional terms to compute before simplifying and applying the final $P\mathcal{L}$:
\begin{align*}
\mathcal{L}P\mathcal{L}Q\mathcal{L}\mathbf{u}_{\mathbf{k}}^0 =&  \mathcal{L}[\mathbf{D}_{\mathbf{k}}(\hat{\mathbf{u}}^0,\tilde{\mathbf{C}}(\hat{\mathbf{u}}^0,\hat{\mathbf{u}}^0))]\\
=&\mathbf{D}_{\mathbf{k}}(\hat{\mathbf{C}}(\mathbf{u}^0,\mathbf{u}^0),\tilde{\mathbf{C}}(\hat{\mathbf{u}}^0,\hat{\mathbf{u}}^0))+\mathbf{D}_{\mathbf{k}}(\hat{\mathbf{u}}^0,\tilde{\mathbf{D}}(\hat{\mathbf{C}}(\mathbf{u}^0,\mathbf{u}^0),\hat{\mathbf{u}}^0))\\
\mathcal{L}\mathcal{L}Q\mathcal{L}\mathbf{u}_{\mathbf{k}}^0 = &\mathcal{L}\mathcal{L}[\mathbf{D}_{\mathbf{k}}(\hat{\mathbf{u}}^0,\tilde{\mathbf{u}}^0)+\mathbf{C}_{\mathbf{k}}(\tilde{\mathbf{u}}^0,\tilde{\mathbf{u}}^0)]\\
=&\mathcal{L}[\mathbf{D}_{\mathbf{k}}(\hat{\mathbf{C}}(\mathbf{u}^0,\mathbf{u}^0),\tilde{\mathbf{u}}^0)+\mathbf{D}_{\mathbf{k}}(\hat{\mathbf{u}}^0,\tilde{\mathbf{C}}(\mathbf{u}^0,\mathbf{u}^0))+\mathbf{D}_{\mathbf{k}}(\tilde{\mathbf{C}}(\mathbf{u}^0,\mathbf{u}^0),\tilde{\mathbf{u}}^0)]\\
=& \mathbf{D}_{\mathbf{k}}(\hat{\mathbf{D}}(\mathbf{C}(\mathbf{u}^0,\mathbf{u}^0),\mathbf{u}^0),\tilde{\mathbf{u}}^0)+2\mathbf{D}_{\mathbf{k}}(\hat{\mathbf{C}}(\mathbf{u}^0,\mathbf{u}^0),\tilde{\mathbf{C}}(\mathbf{u}^0,\mathbf{u}^0))\\
&+ \mathbf{D}_{\mathbf{k}}(\mathbf{u}^0,\tilde{\mathbf{D}}(\mathbf{C}(\mathbf{u}^0,\mathbf{u}^0),\mathbf{u}^0))+\mathbf{D}_{\mathbf{k}}(\tilde{\mathbf{C}}(\mathbf{u}^0,\mathbf{u}^0),\tilde{\mathbf{C}}(\mathbf{u}^0,\mathbf{u}^0)).
\end{align*}Combining these three computed terms with the correct coefficients yields:
\begin{align*}
[3P\mathcal{L}(2P\mathcal{L}-P\mathcal{L})-3\mathcal{L}P\mathcal{L}+\mathcal{L}\mathcal{L}]Q\mathcal{L}\mathbf{u}_{\mathbf{k}}^0 =& \mathbf{D}_{\mathbf{k}}(\hat{\mathbf{u}}^0,\tilde{\mathbf{D}}(-3\hat{\mathbf{C}}(\mathbf{u}^0,\mathbf{u}^0)+3\hat{\mathbf{C}}(\hat{\mathbf{u}}^0,\hat{\mathbf{u}}^0)-3\tilde{\mathbf{C}}(\hat{\mathbf{u}}^0,\hat{\mathbf{u}}^0),\hat{\mathbf{u}}^0))\\
&-6\mathbf{C}_{\mathbf{k}}(\tilde{\mathbf{C}}(\hat{\mathbf{u}}^0,\hat{\mathbf{u}}^0),\tilde{\mathbf{C}}(\hat{\mathbf{u}}^0,\hat{\mathbf{u}}^0))-3\mathbf{D}_{\mathbf{k}}(\hat{\mathbf{C}}(\mathbf{u}^0,\mathbf{u}^0),\tilde{\mathbf{C}}(\hat{\mathbf{u}}^0,\hat{\mathbf{u}}^0))\\
&+\mathbf{D}_{\mathbf{k}}(\hat{\mathbf{D}}(\mathbf{C}(\mathbf{u}^0,\mathbf{u}^0),\mathbf{u}^0),\tilde{\mathbf{u}}^0)+2\mathbf{D}_{\mathbf{k}}(\hat{\mathbf{C}}(\mathbf{u}^0,\mathbf{u}^0),\tilde{\mathbf{C}}(\mathbf{u}^0,\mathbf{u}^0))\\
&+ \mathbf{D}_{\mathbf{k}}(\mathbf{u}^0,\tilde{\mathbf{D}}(\mathbf{C}(\mathbf{u}^0,\mathbf{u}^0),\mathbf{u}^0))+2\mathbf{C}_{\mathbf{k}}(\tilde{\mathbf{C}}(\mathbf{u}^0,\mathbf{u}^0),\tilde{\mathbf{C}}(\mathbf{u}^0,\mathbf{u}^0)).
\end{align*}The $t^3$-term is found once we apply $P\mathcal{L}$ to this expression, yielding:
\begin{align*}
\mathbf{R}_{\mathbf{k}}^3(\hat{\mathbf{u}}) 
=&\mathbf{D}_{\mathbf{k}}(\hat{\mathbf{u}},\tilde{\mathbf{D}}(\hat{\mathbf{u}},\hat{\mathbf{D}}(\hat{\mathbf{u}},\hat{\mathbf{C}}(\hat{\mathbf{u}},\hat{\mathbf{u}})-2\tilde{\mathbf{C}}(\hat{\mathbf{u}},\hat{\mathbf{u}}))+\tilde{\mathbf{D}}(\hat{\mathbf{u}},\tilde{\mathbf{C}}(\hat{\mathbf{u}},\hat{\mathbf{u}})-2\hat{\mathbf{C}}(\hat{\mathbf{u}},\hat{\mathbf{u}})))\\
&\qquad\quad+\tilde{\mathbf{D}}(\tilde{\mathbf{C}}(\hat{\mathbf{u}},\hat{\mathbf{u}}),\tilde{\mathbf{C}}(\hat{\mathbf{u}},\hat{\mathbf{u}})-\hat{\mathbf{C}}(\hat{\mathbf{u}},\hat{\mathbf{u}}))+\tilde{\mathbf{D}}(\hat{\mathbf{C}}(\hat{\mathbf{u}},\hat{\mathbf{u}}),\hat{\mathbf{C}}(\hat{\mathbf{u}},\hat{\mathbf{u}})))\\
&\;\;\;
+3\mathbf{D}_{\mathbf{k}}(\tilde{\mathbf{C}}(\hat{\mathbf{u}},\hat{\mathbf{u}}),\tilde{\mathbf{D}}(\hat{\mathbf{u}},\tilde{\mathbf{C}}(\hat{\mathbf{u}},\hat{\mathbf{u}})-\hat{\mathbf{C}}(\hat{\mathbf{u}},\hat{\mathbf{u}})))
.
\end{align*}

\subsection{Fourth-order term}

We will include up through the fourth-order term in our renormalized reduced order models. The derivation of these models is quite tedious. For this reason, we make use of our symbolic tools described above. The result for the fourth order model is:

\begin{align*}
\mathbf{R}_{\mathbf{k}}^4(\hat{\mathbf{u}})=&e^{t\mathcal{L}}P\mathcal{L}[PLPLPL - 3 PLPLQL - 
 5 PLQLPL  
 +3 PLQLQL  \\&\qquad\qquad\quad\;\;- 
 3 QLPLPL + 
 5 QLPLQL+ 
 3 QLQLPL - QLQLQL]Q\mathcal{L}\mathbf{u}_{\mathbf{k}}^0\\
 =&
\mathbf{D}_{\mathbf{k}}(\hat{\mathbf{u}},\tilde{\mathbf{D}}(\hat{\mathbf{u}},\hat{\mathbf{D}}(\hat{\mathbf{C}}(\hat{\mathbf{u}},\hat{\mathbf{u}}),\hat{\mathbf{C}}(\hat{\mathbf{u}},\hat{\mathbf{u}})-2\tilde{\mathbf{C}}(\hat{\mathbf{u}},\hat{\mathbf{u}}))+3\hat{\mathbf{D}}(\tilde{\mathbf{C}}(\hat{\mathbf{u}},\hat{\mathbf{u}}),\tilde{\mathbf{C}}(\hat{\mathbf{u}},\hat{\mathbf{u}}))\\
&\qquad\qquad\quad+\tilde{\mathbf{D}}(\hat{\mathbf{C}}(\hat{\mathbf{u}},\hat{\mathbf{u}}),2\tilde{\mathbf{C}}(\hat{\mathbf{u}},\hat{\mathbf{u}})-3\hat{\mathbf{C}}(\hat{\mathbf{u}},\hat{\mathbf{u}}))-\tilde{\mathbf{D}}
(\tilde{\mathbf{C}}(\hat{\mathbf{u}},\hat{\mathbf{u}}),\tilde{\mathbf{C}}(\hat{\mathbf{u}},\hat{\mathbf{u}}))\\
&\qquad\qquad\quad+\hat{\mathbf{D}}(\hat{\mathbf{u}},\hat{\mathbf{D}}(\hat{\mathbf{u}},\hat{\mathbf{C}}(\hat{\mathbf{u}},\hat{\mathbf{u}})-3\tilde{\mathbf{C}}(\hat{\mathbf{u}},\hat{\mathbf{u}}))+\tilde{\mathbf{D}}(\hat{\mathbf{u}},3\tilde{\mathbf{C}}(\hat{\mathbf{u}},\hat{\mathbf{u}})-5\hat{\mathbf{C}}(\hat{\mathbf{u}},\hat{\mathbf{u}})))\\
&\qquad\qquad\quad+\tilde{\mathbf{D}}(\hat{\mathbf{u}},\hat{\mathbf{D}}(\hat{\mathbf{u}},5\tilde{\mathbf{C}}(\hat{\mathbf{u}},\hat{\mathbf{u}})-3\hat{\mathbf{C}}(\hat{\mathbf{u}},\hat{\mathbf{u}}))+\tilde{\mathbf{D}}(\hat{\mathbf{u}},3\hat{\mathbf{C}}(\hat{\mathbf{u}},\hat{\mathbf{u}})-\tilde{\mathbf{C}}(\hat{\mathbf{u}},\hat{\mathbf{u}}))))\\
&\qquad\quad+\tilde{\mathbf{D}}(\hat{\mathbf{C}}(\hat{\mathbf{u}},\hat{\mathbf{u}}),\hat{\mathbf{D}}(\hat{\mathbf{u}},3\hat{\mathbf{C}}(\hat{\mathbf{u}},\hat{\mathbf{u}})-5\tilde{\mathbf{C}}(\hat{\mathbf{u}},\hat{\mathbf{u}}))+\tilde{\mathbf{D}}(\hat{\mathbf{u}},\tilde{\mathbf{C}}(\hat{\mathbf{u}},\hat{\mathbf{u}})-3\hat{\mathbf{C}}(\hat{\mathbf{u}},\hat{\mathbf{u}})))\\
&\qquad\quad+\tilde{\mathbf{D}}(\tilde{\mathbf{C}}(\hat{\mathbf{u}},\hat{\mathbf{u}}),\hat{\mathbf{D}}(\hat{\mathbf{u}},3\tilde{\mathbf{C}}(\hat{\mathbf{u}},\hat{\mathbf{u}})-\hat{\mathbf{C}}(\hat{\mathbf{u}},\hat{\mathbf{u}}))+\tilde{\mathbf{D}}(\hat{\mathbf{u}},5\hat{\mathbf{C}}(\hat{\mathbf{u}},\hat{\mathbf{u}})-3\tilde{\mathbf{C}}(\hat{\mathbf{u}},\hat{\mathbf{u}}))))\\
&-4\mathbf{D}_{\mathbf{k}}(\tilde{\mathbf{C}}(\hat{\mathbf{u}},\hat{\mathbf{u}}),\tilde{\mathbf{D}}(\hat{\mathbf{C}}(\hat{\mathbf{u}},\hat{\mathbf{u}}),\hat{\mathbf{C}}(\hat{\mathbf{u}},\hat{\mathbf{u}})-\tilde{\mathbf{C}}(\hat{\mathbf{u}},\hat{\mathbf{u}}))+\tilde{\mathbf{D}}(\tilde{\mathbf{C}}(\hat{\mathbf{u}},\hat{\mathbf{u}}),\tilde{\mathbf{C}}(\hat{\mathbf{u}},\hat{\mathbf{u}}))\\
&\qquad\qquad\qquad\quad+\tilde{\mathbf{D}}(\hat{\mathbf{u}},\hat{\mathbf{D}}(\hat{\mathbf{u}},\hat{\mathbf{C}}(\hat{\mathbf{u}},\hat{\mathbf{u}})-2\tilde{\mathbf{C}}(\hat{\mathbf{u}},\hat{\mathbf{u}}))+\tilde{\mathbf{D}}(\hat{\mathbf{u}},\tilde{\mathbf{C}}(\hat{\mathbf{u}},\hat{\mathbf{u}})-2\hat{\mathbf{C}}(\hat{\mathbf{u}},\hat{\mathbf{u}}))))\\
&-3\mathbf{D}_{\mathbf{k}}(\tilde{\mathbf{D}}(\hat{\mathbf{u}},\hat{\mathbf{C}}(\hat{\mathbf{u}},\hat{\mathbf{u}})),\tilde{\mathbf{D}}(\hat{\mathbf{u}},\hat{\mathbf{C}}(\hat{\mathbf{u}},\hat{\mathbf{u}})-2\tilde{\mathbf{C}}(\hat{\mathbf{u}},\hat{\mathbf{u}})))\\
&-3\mathbf{D}_{\mathbf{k}}(\tilde{\mathbf{D}}(\hat{\mathbf{u}},\tilde{\mathbf{C}}(\hat{\mathbf{u}},\hat{\mathbf{u}})),\tilde{\mathbf{D}}(\hat{\mathbf{u}},\tilde{\mathbf{C}}(\hat{\mathbf{u}},\hat{\mathbf{u}})))
.
\end{align*}

\end{document}